\documentclass[opre,sglanonrev]{informs4}

\OneAndAHalfSpacedXI

\usepackage{algorithm}

\usepackage{amsmath,amssymb,amsfonts}
\usepackage{amstext}
\usepackage{multirow}
\usepackage{algorithmic}
\usepackage{amsopn}
\usepackage{threeparttable}
\usepackage{hyperref}
\usepackage{subfigure}
\usepackage{lineno}

\usepackage{rotating}
\usepackage{fancyvrb}

\usepackage[sort&compress]{natbib}
\bibpunct[, ]{[}{]}{,}{n}{}{,}%
\def\bibfont{\small}%
\EquationsNumberedThrough    

\TheoremsNumberedThrough     
\ECRepeatTheorems  %

\MANUSCRIPTNO{}

\begin{document}

 \RUNAUTHOR{Zhang et al.} 

\RUNTITLE{Iterative Minimax Games with Coupled Linear Constraints}

\TITLE{Iterative Minimax Games with Coupled Linear Constraints
	\thanks{This research is supported by National Natural Science Foundation of China under the Grants (Nos. 12471294, 12021001, 11991021), and by the Strategic Priority Research Program of Chinese Academy of Sciences (No. XDA27000000).}}

\ARTICLEAUTHORS{%
\AUTHOR{Huiling Zhang}

\AFF{Department of Mathematics, Shanghai University, Shanghai 200444, People’s Republic of China; 
	LSEC, ICMSEC, Academy of Mathematics and Systems Science,
	Chinese Academy of Sciences, Beijing 100190, China, \EMAIL{zhanghl1209@shu.edu.cn}}

\AUTHOR{Zi Xu\thanks{Corresponding author.}}

\AFF{Department of Mathematics, Shanghai University, Shanghai 200444, People’s Republic of China; \break Newtouch Center for Mathematics of Shanghai University, Shanghai 200444, People’s Republic of China, \EMAIL{xuzi@shu.edu.cn}.}

\AUTHOR{Yu-Hong Dai}

\AFF{LSEC, ICMSEC, Academy of Mathematics and Systems Science,
	Chinese Academy of Sciences, Beijing 100190, China, \EMAIL{dyh@lsec.cc.ac.cn}}

%
}

\ABSTRACT{%
	The study of nonconvex minimax games has gained significant momentum in machine learning and decision science communities due to their fundamental connections to adversarial training scenarios. This work develops a primal-dual alternating proximal gradient (PDAPG) algorithm framework for resolving iterative minimax games featuring nonsmooth nonconvex objectives subject to coupled linear constraints. We establish rigorous convergence guarantees for both nonconvex-strongly concave and nonconvex-concave game configurations, demonstrating that PDAPG achieves an $\varepsilon$-stationary solution within $\mathcal{O}\left( \varepsilon ^{-2} \right)$ iterations for strongly concave settings and $\mathcal{O}\left( \varepsilon ^{-4} \right)$ iterations for concave scenarios. Our analysis provides the first known iteration complexity bounds for this class of constrained minimax games, particularly addressing the critical challenge of coupled linear constraints that induce inherent interdependencies among strategy variables. The proposed game-theoretic framework advances existing solution methodologies by simultaneously handling nonsmooth components and coordinated constraint structures through alternating primal-dual updates.}
%
%


\maketitle

\section{Introduction}\label{sec1}
In recent times, iterative minimax games with coupled linear constraints have emerged as a prominent area of research within multiple fields such as machine learning, signal processing, statistics, and economics \cite{{Berger,Facchinei,Ghosh,Sanjabi,Wu,Yang}}. A multitude of practical applications can be cast into the form of minimax problems. For instance, in signal processing, the power control and transceiver design problem \cite{Lu}, reinforcement learning \cite{Cai,Moriarty,Qiu}, adversarial attacks on deep neural networks (DNNs) \cite{Chen}, distributed nonconvex optimization \cite{Giannakis,Mateos}, and robust learning across multiple domains \cite{Qian} can all be formulated as: 
\begin{equation}\label{min} 
\min \limits_{x\in \mathcal{X}}\max \limits_{y\in \mathcal{Y}}f(x,y), 
\end{equation} 
Here, $\mathcal{X}\subseteq \mathbb{R}^n$ and $\mathcal{Y}\subseteq \mathbb{R}^m$ represent nonempty convex and compact sets. The function $f(x,y): \mathcal{X} \times \mathcal{Y} \rightarrow \mathbb{R}$ is smooth but may be nonconvex with respect to $x$. We term \eqref{min} a (non)convex-(non)concave minimax problem based on whether $f(x, y)$ is (non)convex in $x$ and (non)concave in $y$. Numerous studies have been carried out on different types of minimax problems, including (strongly) convex-(strongly) concave \cite{{Lan2016,Mokhtari,Nemirovski,Nesterov2007,Thek2022}}, nonconvex-(strongly) concave \cite{{Lin2019,Lin2020,Lu,Pan,Shen,Xu,Xu21,Xu2021,xu22zeroth,Zhang}}, and nonconvex-nonconcave \cite{{Bohm,DiakonikolasEG,Doan,Lee,Sanjabi18,SongOptDE}} cases. Zhang et al. \cite{{Zhang2021}} put forward a general acceleration framework for solving the nonconvex-strongly concave minimax problem \eqref{min}. This algorithm can identify an $\varepsilon$-stationary point within $\mathcal{O}\left(\sqrt{\kappa_y}\varepsilon ^{-2} \right)$, which stands as the optimal complexity bound under the given setting \cite{{Li21}}.

While minimax optimization has been extensively studied in unconstrained settings, practical scenarios with coupled linear constraints introduce fundamentally new challenges in equilibrium-seeking dynamics. A prominent example arises in adversarial network flow games \cite{Tsaknakis}, where strategic agents interact under shared resource limitations. Specifically, consider a directed graph $G=(V, E)$, where regular users aim to route flow $x\in \mathcal{X}$ from source  $s$ to sink $t$ with minimal cost, while an adversary injects perturbative flow $y\in \mathcal{Y}$ to maximize regular users' expenditures. This adversarial interaction is governed by the bilinear coupling constraint $\mathbf{x}+\mathbf{y} \leq \mathbf{p}$  forming the minimax game:
\begin{align}
	\max _{y\in \mathcal{Y},z\in\mathbb{R}_{+}^{|E|}}  \min _{x\in \mathcal{X}} & \quad \sum_{(i, j) \in E} q_{i j}\left(x_{i j}+y_{i j}\right) x_{i j} -\frac{\eta}{2}\|y\|^2\nonumber\\
	\mbox { s.t. } & \quad \mathbf{x}+\mathbf{y}+\mathbf{z}  =\mathbf{p}.\label{testprobn1}
\end{align}
where $\mathcal{X}= \{x\mid 0 \leq \mathbf{x} \leq \mathbf{p}, \sum_{(i, t) \in E} x_{i t}=r_t, \sum_{(i, j) \in E} x_{i j}-\sum_{(j, k) \in E} x_{j k}=0, \forall j \in V \backslash\{s, t\}\}$, $\mathcal{Y}=\{y\mid 0 \leq \mathbf{y} \leq \mathbf{p}, \sum_{(i, j) \in E} y_{i j}=b\}$ which enforces flow conservation and sink demand $r_t$,  $V$ and $E$ is the set of vertices and edges of a directed graph $G=(V, E)$, respectively, $s$ and $t$ denote the source and sink node respectively, and $r_t$ is the total flow to sink node $t$, $\mathbf{x}=\left\{x_e\right\}_{e \in E} \in \mathbb{R}^{|E|}$ and $\mathbf{p}=\left\{p_e\right\}_{e \in E} \in \mathbb{R}^{|E|}$ denote the vectors of flows and the capacities at edge $e$, respectively, and $\mathbf{y}=\left\{y_e\right\}_{e \in E}$  denote the set of flows controlled by the adversary, $\eta>0$, $q_{ij}>0$ is the unit cost at edge $(i,j)\in E$ and $b > 0$ is its total budget.  $\mathcal{Y}$ restricts adversarial budgets. The quadratic penalty $-\frac{\eta}{2}\|y\|^2$ regularizes adversary's intervention intensity.

%

Recent advances in distributionally robust optimization (DRO) reveal critical demands for iterative minimax games with coupled linear constraints beyond convex settings. The standard DRO framework \cite{Namkoong,Zhang24} minimizes the worst-case empirical risk via
\begin{align*}
	\min _{x \in \mathbb{R}^d} \max _{y \in Y} \frac{1}{n} \sum_{i=1}^n y_i \ell_i(x), 
\end{align*}
where nonconvex losses $\ell_i(x)$  induce complex adversarial dynamics. Practical implementations further require coupled linear constraints $Ax+By=c$ to model strategic interactions between agents under shared resources. For instance, in portfolio optimization, investment allocations $x$ and market scenario weights $y$ must satisfy budget conservation across asset classes; in adversarial image recognition, feature selection $x$ and class reweighting $y$ are linearly correlated to balance model robustness and fairness. These scenarios motivate the constrained minimax formulation:
\begin{align}
	\min _{x \in \mathbb{R}^d} \max _{y \in Y}&\quad  \frac{1}{n} \sum_{i=1}^n y_i \ell_i(x),\nonumber\\
	\mbox { s.t. } & \quad Ax+By=c.\label{test-dro}
\end{align}
where the linear coupling $Ax+By=c$ explicitly intertwines the agent's decision $x$ with the adversary's strategy $y$, transforming conventional DRO into an equilibrium-seeking game. This equilibrium structure distinguishes the problem from conventional unconstrained minimax settings, demanding novel algorithmic approaches to resolve the tension between strategic competition and cooperative constraint satisfaction.

In this paper, we consider the following general nonsmooth nonconvex minimax game problems with coupled linear constraints:
\begin{equation}\label{problem:1}
	\min_{x\in\mathcal{X}}\max\limits_{\substack{ y\in\mathcal{Y}\\Ax+By=c}}\{f(x,y)+h(x)-g(y)\},
\end{equation}
where $\mathcal{X}\subseteq \mathbb{R}^n$ and $\mathcal{Y}\subseteq \mathbb{R}^m$ are nonempty convex and compact sets, $f(x,y): \mathcal{X} \times \mathcal{Y} \rightarrow \mathbb{R}$ is a smooth function and possibly nonconvex in $x$, (strongly) concave in $y$, $h(x):\mathbb{R}^{n}\rightarrow \mathbb{R}$ and $g(y): \mathbb{R}^{m}\rightarrow \mathbb{R}$ are some proper convex continuous possibly nonsmooth functions, $A\in\mathbb{R}^{p\times n}$, $B\in\mathbb{R}^{p\times m}$ and $c\in\mathbb{R}^{p}$. We use $[Ax+By]_i= [c]_i$ to denote the $i$th constraint, and define $\mathcal{K}:= \{1,\dots , p\}$ as the index set of the constraints. Problem \eqref{problem:1} is more challenging than \eqref{min}, and it is NP-hard to find its globally optimal solutions even when $f(x, y)$ is strongly-convex in $x$ and strongly-concave in $y$~\cite{Tsaknakis}. 

Existing methods for linearly constrained minimax games predominantly adopt multi-loop architectures.
In \cite{Tsaknakis}, a family of efficient algorithms named multiplier gradient descent (MGD) are proposed for solving smooth minimax problem with coupled linear inequality constraints, which achieves the iteration complexity of $\tilde{\mathcal{O}}\left( \varepsilon ^{-2} \right)$ when $f(x,y)$ is strongly convex with respect to $x$ and strongly concave with respect to $y$.
Dai et al. \cite{Dai} proposed a proximal gradient multi-step ascent descent method (PGmsAD) for nonsmooth convex-concave bilinear coupled minimax problems with linearly equality constraints when $\mathcal{X}=\mathbb{R}^n$ and $\mathcal{Y}=\mathbb{R}^m$, which achieves the iteration complexity of $\mathcal{O}(\varepsilon^{-2}\log(1/\varepsilon))$.

Multi-loop algorithms for minimax optimization typically involve nested iterations, where each outer update requires solving an inner subproblem to high precision. For instance, classical approaches to constrained minimax problems often employ sequential linear approximations of nonlinear functions, solving linearized subproblems iteratively. This nested structure introduces significant computational overhead, particularly in large-scale networks where dimensionality exacerbates the cost of subproblem resolution. Similarly, in dynamic game-theoretic frameworks, policy iteration algorithms for solving nonzero-sum games with coupled constraints often alternate between value function evaluation and policy improvement steps, requiring multiple inner iterations to stabilize intermediate solutions. These methods, though theoretically sound, face practical bottlenecks in scenarios demanding real-time adaptability or deployment in resource-constrained environments. These limitations motivate the development of single-loop algorithms that unroll nested iterations into a unified update framework. By eliminating inner-loop computations and leveraging closed-form operations (e.g., gradient-based adjustments with adaptive step sizes), single-loop methods can achieve comparable convergence guarantees while drastically reducing per-iteration complexity—a critical advantage for applications in distributed systems, real-time control, or high-dimensional optimization.

Notably, single-loop algorithms that update variables via closed-form operations remain absent for convex-(strongly) concave cases with coupled constraints, despite their advantages in iteration efficiency and implementation simplicity \cite{Zhang}. This gap motivates our development of a provably convergent single-loop framework for iterated minimax games with linear coupling. To the best of our knowledge, there is no single-loop algorithm with a theoretical complexity guarantee for nonconvex-(strongly) concave minimax problems with coupled linear constraints up to now.

\subsection{Motivated Applications}\label{app}
We give the following problems which can be formulated by \eqref{problem:1} as motivated examples for this paper.

\textbf{Multi-level programming problem.} Multi-level programming problem has been a focus of research in supply chain management, traffic and transportation network design, energy management, and many other research fields. It inherently embodies sequential minimax dynamics with coupled linear constraints. Consider a hierarchy of  $n$ strategic agents where each agent $i$ controls variables $x^i\in X^i\subset \mathbb{R}^{n_i}$ to maximize their objective $f_i\left(x^1, x^2, \ldots, x^n\right) $ subject to interdependencies induced by preceding agents' decisions. 

Specifically, agent $1$ first selects  $x^1\in X_1$, followed by agent 2 choosing $x^2\in X^2$ conditioned on $x_1$, and so forth—a process mirroring iterative adversarial updates in minimax games. When $n=2$ with diametrically opposed objectives ($f_2\left(x^1, x^2\right) = -f_1\left(x^1, x^2\right) $) and linear constraint functions $g^1, g^2$, this hierarchy reduces to a linear coupled constrained minimax problem~\cite{Falk}:
\begin{align*}
	\max_{\substack{x^1 \in X^1\\g^1\left(x^1\right) \geq 0}}\min\limits_{\substack{ x^2 \in X^2\\g^2\left(x^1, x^2\right) \geq 0}}f_1\left(x^1, x^2\right) .
\end{align*}

\textbf{Generalized absolute value equations (GAVEs).} The GAVE is a popular nonsmooth NP-hard problem which can be formulated as follows:
\begin{align}\label{gave}
A x+B|x|=b,
\end{align}
where $A, B \in \mathbb{R}^{m \times n}, b \in \mathbb{R}^m$, and for $x=\left(x_1, \ldots, x_n\right) \in \mathbb{R}^n,|x|:=\left(\left|x_1\right|, \ldots,\left|x_n\right|\right) \in$ $\mathbb{R}_{+}^n$. As an important tool in optimization, GAVE is broadly applied to address problems in multiple fields, such as nonnegative constrained least squares problems, quadratic programming, and bimatrix games \cite{Dai}. In \cite{Dai}, it is proved that GAVE is equivalent to the following convex-concave minimax problem with linearly equality constraints:
\begin{align}
	\min _{x^+ \in \mathbb{R}_{+}^n} \max _{z \in \mathbb{R}_{+}^n, y \in \mathbb{R}^m} & (b-(A+B) x^+)^T y \nonumber\\
	\text { subject to } & x^+-(B-A)^T y-z=0, \label{testprob2}
\end{align}
where $x^+=(\max\{0,x_1\}, \ldots, \max\{0,x_n\})$.
Clearly, problem \eqref{testprob2} is a special case of \eqref{problem:1}.

\textbf{Robust portfolio optimization.} Modern portfolio optimization exemplifies iterative minimax dynamics under coupled resource constraints, where investors strategically balance returns and risks under shared budgetary limitations. Consider an investor allocating wealth proportions $x \in \mathcal{X} = [0, 1]^n$  across $n$ assets to maximize expected return $\mu^\top x$ while mitigating tail risks, subject to the linear budgetary $e^\top x=1$. Traditional mean-variance frameworks \cite{Leung20} optimize the Sharpe ratio $ -\frac{\mu^\top x-r_f}{\sqrt{x^\top V x }}$ but suffer from symmetric risk treatment, ignoring asymmetric loss aversion under non-normal market conditions. To address this, \cite{Wang} reformulates the problem as a constrained minimax game via CVaR, introducing adversarial dynamics between return-seeking and risk mitigation:
\begin{align*}
	\min _{x\in\mathcal{X}, \rho} \max _\lambda\quad &\lambda\left(\rho+\frac{1}{N(1-\theta)} \sum_{i=1}^N\left(- \xi_j^T x-\rho\right)^{+}\right)-\left(\mu^T x-r_f\right) \\
	\text { s.t. } &e_n^T x=1\\
	&\lambda\left(\rho+\frac{1}{N(1-\theta)} \sum_{i=1}^N\left(- \xi_j^T x-\rho\right)^{+}\right)=\mu^T x-r_f .
\end{align*}
where $\lambda > 0$ modulates risk sensitivity, $\rho$ calibrates loss thresholds via CVaR approximation, and $\xi_j$ represent observed returns. Here, the investor $(x, \rho)$ minimizes combined costs—enhancing returns while suppressing CVaR—while the risk modeler ($\lambda$) maximizes penalties to expose portfolio fragility, all under two coupled constraints: the budgetary coupling $e^\top x=1$ ensures resource coherence, while the risk-return parity equality forces CVaR penalties to scale with realized excess returns. Each iteration adversarially updates $x$ (rebalancing assets to exploit high-$\mu$ opportunities), $\rho$ (tightening loss thresholds to identify vulnerabilities), and $\lambda$ (calibrating penalty intensity to balance risk-return tradeoffs). This creates a Nash equilibrium where no agent can unilaterally improve without violating shared constraints—a hallmark of nonconvex minimax games with linear interdependencies. The
$\max(\cdot)$ operator in CVaR and equality constraint jointly induce nonconvexity, precluding standard gradient-based methods and necessitating iterative algorithms that resolve competitive objectives alongside cooperative constraint satisfaction.

\subsection{Contributions}
This work develops a single-loop algorithmic framework for nonsmooth nonconvex minimax optimization under coupled linear equality constraints, advancing both theoretical and computational frontiers. Our primary contributions are fourfold. First, we establish strong duality with respect to $y$ for nonconvex-concave problems under feasible constraint qualifications—a novel result for linearly constrained minimax games, previously unexplored beyond convex settings. Leveraging this duality, we propose the Primal-Dual Alternating Proximal Gradient (PDAPG) algorithm, a single-loop method that alternately updates primal ($x$) and dual $(y,\lambda )$ variables via proximal steps, circumventing nested iterations. Theoretically, PDAPG achieves an $\mathcal{O}\left( \varepsilon ^{-2} \right)$ (resp. $\mathcal{O}\left( \varepsilon ^{-4} \right)$) iteration complexity for nonconvex-strongly concave (resp. nonconvex-concave) problems with coupled constraints, matching the lower-bound $\mathcal{O}\left( \varepsilon ^{-2} \right)$ for unconstrained nonconvex-strongly concave cases and aligning with state-of-the-art rates of AGP \cite{Xu} and smoothed GDA \cite{Zhang} for unconstrained nonconvex-concave settings. Crucially, PDAPG attains these guarantees without sacrificing the single-loop structure, even when resolving linear coupling constraints—a marked improvement over multi-loop methods like MGD \cite{Tsaknakis} and PGmsAD \cite{Dai}. This work thus bridges a critical gap in constrained minimax optimization, delivering the first provably efficient single-loop solver for coupled nonconvex-(strongly) concave games while preserving computational tractability.

\subsection{Related Works}
In this subsection, we first give a brief review of the related methods for solving minimax optimization problem \eqref{min}. 

Under the convex-concave setting of \eqref{min}, there are many existing works. For instance, Nemirovski \cite{Nemirovski} proposed a mirror-prox method with complexity of $\mathcal{O}(\varepsilon^{-1})$ in terms of duality gap when $\mathcal{X}$ and $\mathcal{Y}$ are bounded sets. Nesterov \cite{Nesterov2007} proposed a dual extrapolation algorithm with complexity of $\mathcal{O}(\varepsilon^{-1})$. Mokhtari et al. \cite{Mokhtari} established an overall complexity of $\mathcal{O}(\kappa\log(1/\varepsilon))$ for both extra-gradient (EG) and optimistic gradient descent ascent (OGDA) methods in bilinear and strongly convex-strongly concave settings. Thekumparampil et al. \cite{Thek2022} proposed a lifted primal-dual method for  bilinearly coupled strongly convex-strongly concave minimax problem, which achieves the iteration complexity of $\mathcal{O}(\sqrt{\kappa}\log(1/\varepsilon))$. For more related results, we refer to \cite{Chen2017,Chen14,Dang,Gidel2019,He2016,Lan2016,Lin2020,Ouyang2015,Ouyang2019,Tominin,Zhang22} and the references therein.

For nonconvex-strongly concave minimax problem \eqref{min}, various algorithms have been proposed in recent works \cite{Bot,Jin,Lin2019,Lin2020,Lu,Rafique}, and all of them can achieve the iteration complexity of $\tilde{\mathcal{O}}\left( \kappa_y^2\varepsilon ^{-2} \right)$ in terms of stationary point of $\Phi (\cdot) = \max_{y\in \mathcal{Y}} f(\cdot, y)$ (when $\mathcal{X}=\mathbb{R}^n$), or stationary point of $f(x,y)$,  where $\kappa_y$ is the condition number for $f(x,\cdot)$. Furthermore, Zhang et al. \cite{Zhang2021} proposed a generic acceleration framework which can improve the iteration complexity to $\mathcal{O}\left(\sqrt{\kappa_y}\varepsilon ^{-2} \right)$.

For general nonconvex-concave minimax problem \eqref{min}, there are two types of algorithms, i.e., multi-loop algorithms and single-loop algorithms. Various multi-loop algorithms have been proposed in \cite{Kong,Lin2020,Nouiehed,Ostro,Rafique,Thek2019}. The best known iteration complexity in terms of stationary point of $f(x,y)$ for multi-loop algorithms is $\tilde{\mathcal{O}}\left( \varepsilon ^{-2.5} \right)$, which can be achieved by \cite{Lin2020}. All these multi-loop algorithms either utilize multiple accelerated gradient ascent steps to update $y$ that can solve the inner subproblem exactly or inexactly, or further do similar acceleration for $x$'s update by adding regularization terms to the inner objective function, and thus are relatively complicated to be implemented.
On the other hand, fewer studies focus on single-loop algorithms for nonconvex-concave minimax problems. One classical method is the gradient descent-ascent (GDA) method, which employs a gradient descent step to update $x$ and a gradient ascent step to update $y$ simultaneously at each iteration. Many other GDA-type algorithms have been analyzed in \cite{Chambolle,Daskalakis17,Daskalakis18,Gidel2019,Gidel18-2,Ho2016}. Lin et al. \cite{Lin2019} proved that the iteration complexity of GDA  to find an $\varepsilon$-stationary point of $\Phi (\cdot) = \max_{y\in \mathcal{Y}} f(\cdot, y)$ is $\tilde{\mathcal{O}} (\varepsilon^{-6})$ when $\mathcal{X}=\mathbb{R}^n$ and $\mathcal{Y}$ is a convex compact set. Lu et al. \cite{Lu} proposed the hybrid block successive approximation (HiBSA) algorithm, which can obtain the iteration complexity of $\tilde{\mathcal{O}}\left( \varepsilon ^{-4} \right)$ in terms of $\varepsilon$-stationary point of $f(x,y)$ when $\mathcal{X}$ and $\mathcal{Y}$ are convex compact sets. However, to update $y$, HiBSA has to solve the subproblem of maximizing the original function $f(x, y)$ plus some regularization terms for a fixed $x$. Pan et al. \cite{Pan} proposed a new alternating gradient projection algorithm for nonconvex-linear minimax problem with the iteration complexity of $\mathcal{O}\left( \varepsilon ^{-3} \right)$. Xu et al. \cite{Xu} proposed a unified single-loop alternating gradient projection (AGP) algorithm for solving nonconvex-concave and convex-nonconcave minimax problems, which can find an $\varepsilon$-stationary point with the iteration complexity of $\mathcal{O}\left( \varepsilon ^{-4} \right)$ for smooth and nonsmooth cases. Zhang et al. \cite{Zhang} proposed a smoothed GDA algorithm which also achieves the iteration complexity of $\mathcal{O}\left( \varepsilon ^{-4} \right)$ for general nonconvex-concave minimax problems, and can achieve the iteration complexity of $\mathcal{O}\left( \varepsilon ^{-2} \right)$ for a class of special nonconvex-linear minimax problems.

For nonconvex-nonconcave minimax problem \eqref{min}, Sanjabi et al. \cite{Sanjabi18} proposed a multi-step gradient descent-ascent algorithm which finds an $\varepsilon$-first order Nash equilibrium in $\mathcal{O}(\varepsilon^{-2}\log\varepsilon^{-1})$ iterations when a one-sided Polyak-{\L}ojasiewicz (PL) condition is satisfied. Yang et al. \cite{Yang} showed the alternating GDA algorithm converges globally at a linear rate for a subclass of nonconvex-nonconcave objective functions satisfying a two-sided PL inequality. Song et al. \cite{SongOptDE} proposed an optimistic dual extrapolation (OptDE) method with the iteration complexity of $\mathcal{O}\left( \varepsilon^{-2} \right)$, if a weak solution exists. Hajizadeh et al. \cite{Hajizadeh} showed that a damped version of extra-gradient method linearly converges to an $\varepsilon$-stationary point of nonconvex-nonconcave minimax problems that the nonnegative interaction dominance condition is satisfied. 
For more related results, we refer to \cite{Bohm,Chinchilla,Cai,DiakonikolasEG,Doan,Fiez,Grimmer,Jiang,Lee,Ostrovskii}.

\subsection{Organization}
The rest of this paper is organized as follows. In Section 2, we first establish the strong duality with respect to $y$ under some feasibility assumption for nonconvex-concave minimax problem \eqref{problem:1} with linearly coupled equality constraints. Then, we propose a  primal-dual alternating proximal gradient (PDAPG) algorithm for nonsmooth nonconvex-(strongly) concave minimax problem with coupled linear constraints, and then prove its iteration complexity. Numerical results in Section 3 show the efficiency of the proposed algorithm. Some conclusions are made in the last section.
\subsection{Notations}
For vectors, we use $\|\cdot\|$ to represent the Euclidean norm and its induced matrix norm; $\left\langle x,y \right\rangle$ denotes the inner product of two vectors of $x$ and $y$. $[x]_i$ denotes the $i$th component of vector $x$. Let $\mathbb{R}^p$ denote the Euclidean space of dimension $p$ and $\mathbb{R}_{+}^p$ denote the nonnegative orthant in $\mathbb{R}^p$. We use $\nabla_{x} f(x,y,z)$ (or $\nabla_{y} f(x, y,z)$, $\nabla_{z} f(x, y,z)$) to denote the partial derivative of $f(x,y,z)$ with respect to $x$ (or $y$, $z$) at point $(x, y,z)$, respectively.
We use the notation $\mathcal{O} (\cdot)$ to hide only absolute constants which do not depend on any problem parameter, and $\tilde{\mathcal{O}}(\cdot)$ notation to hide only absolute constants and log factors. A continuously differentiable function $f(\cdot)$ is called $L$-smooth if there exists a constant $L> 0$ such
that for any given $x,y \in \mathcal{X}$, we have
$
\|\nabla f(x)-\nabla f(y)\|\le L\|x-y\|.
$
A continuously differentiable function $f(\cdot)$ is called $\mu$-strongly concave if there exists a constant $\mu> 0$ such
that for any $x,y \in \mathcal{X}$, we have
$$
f(y)\le f(x)+\langle \nabla f(x),y-x \rangle-\frac{\mu}{2}\|y-x\|^2.
$$

\section{Nonconvex-(Strongly) Concave Minimax Problems with Coupled Linear Constraints}\label{sec2}

\subsection{Dual Problem and Strong Duality for $y$}
By using the Lagrangian function of problem \eqref{problem:1}, we obtain the dual problem of \eqref{problem:1}, i.e.,
\begin{align}\label{problem:2}
	\min_{\lambda} \min_{x\in\mathcal{X}}\max_{y\in\mathcal{Y}}\left\{\mathcal{L}(x,y,\lambda)+h(x)-g(y)\right\},\tag{D}
\end{align}
where $\mathcal{L}(x,y,\lambda):=f(x,y)-\lambda^\top(Ax+By-c)$. 
The strong duality of problem \eqref{problem:1} with respect to $y$ is established in the following theorem. For the sake of completeness, we give its proof, although it is similar to that of Theorem 3.2 in \cite{Tsaknakis}.

\begin{theorem}\label{dual}
	Suppose $f(x,y)$ is a concave function with respect to $y$, $\mathcal{Y}$ is a convex and compact set. Then the strong duality of problem \eqref{problem:1} with respect to $y$ holds, i.e.,
	\begin{align*}
		&\min_{x\in\mathcal{X}}\max_{\substack{y\in\mathcal{Y}\\Ax+By= c}}\{f(x,y)+h(x)-g(y)\}
		=\min_{\lambda}\min_{x\in\mathcal{X}}\max_{y\in\mathcal{Y}}\left\{\mathcal{L}(x,y,\lambda)+h(x)-g(y)\right\},
	\end{align*}
	and a point $(x^*,y^*)$ is a stationary point of problem \eqref{problem:1} if and only if there exists a $\lambda^*$ such that $(x^*,y^*,\lambda^*)$ is a stationary point of problem \eqref{problem:2}.
\end{theorem}

\begin{proof}{Proof}
	
	By Proposition 5.3.1 in \cite{Bertsekas}, we first obtain
	\begin{align}\label{dual4}
		\max_{y\in\mathcal{Y}}\min_{\lambda}\left\{\mathcal{L}(x,y,\lambda)+h(x)-g(y)\right\}
		=&\max_{\substack{y\in\mathcal{Y}\\Ax+By=  c}}\left\{f(x,y)+h(x)-g(y)\right\},
	\end{align}
	which further implies that 
	\begin{align}
		\min_{x\in\mathcal{X}}\max\limits_{\substack{ y\in\mathcal{Y}\\Ax+By= c}}\{f(x,y)+h(x)-g(y)\}=\min\limits_{x\in\mathcal{X}}\max\limits_{y\in\mathcal{Y}}\min\limits_{\lambda}\left\{\mathcal{L}(x,y,\lambda)+h(x)-g(y)\right\}.\label{thm2.2:1}
	\end{align}
	Moreover, since $f(x,y)$ is a concave function with respect to $y$, by Sion’s minimax theorem (Corollary 3.3 in \cite{Sion}), $\forall x\in\mathcal{X}$, we have
	\begin{align}\label{dual:5}
		&\max_{y\in\mathcal{Y}}\min_{\lambda}\left\{\mathcal{L}(x,y,\lambda)+h(x)-g(y)\right\}
		=\min_{\lambda}\max_{y\in\mathcal{Y}}\left\{\mathcal{L}(x,y,\lambda)+h(x)-g(y)\right\}.
	\end{align}
	The proof is then completed by combining \eqref{dual4} and \eqref{dual:5}.
	%
\end{proof}

Note that a previous strong duality result is proved in Theorem 3.4 in \cite{Tsaknakis} when a special case of \eqref{problem:1} is studied, i.e., $f(x,y)$ is strongly convex with respect to $x$ and strongly concave with respect to $y$, $h(x)=0$ and $g(y)=0$. Theorem \ref{dual} establishes the strong duality for a class of more general nonsmooth nonconvex-concave minimax problems.

For any given proper convex function $h(\cdot)$ and a convex compact set $\mathcal{Z}$, we define the proximity operator as 
$
\operatorname{Prox}_{h,\mathcal{Z}}^{\alpha}(\upsilon):=\arg\min\limits_{z\in \mathcal{Z}} h(z)+\frac{\alpha}{2}\|z-\upsilon \|^2,
$
where $\alpha>0$.  We then define the $\varepsilon$-stationary point of problem \eqref{problem:1} as follows.
\begin{definition}\label{gap-f}
	For some given $\alpha>0$, $\beta>0$ and $\gamma>0$, denote 	\begin{equation*}
		\nabla G^{\alpha,\beta,\gamma}\left( x,y,\lambda \right) :=\left( \begin{array}{c}
			\alpha\left( x-\operatorname{Prox}_{h,\mathcal{X}}^{\alpha}\left( x-\frac{1}{\alpha}\nabla _x\mathcal{L}\left( x,y,\lambda  \right) \right) \right)\\
			\beta\left( y-\operatorname{Prox}_{g,\mathcal{Y}}^{\beta}\left( y+ \frac{1}{\beta} \nabla _y\mathcal{L}\left( x,y,\lambda \right) \right) \right)\\
			\nabla_{\lambda}\mathcal{L}(x,y,\lambda)\\
		\end{array} \right) .
	\end{equation*}
	Then for any given $\varepsilon>0$, if there exists a $\lambda$ such that $\nabla G^{\alpha,\beta,\gamma}\left( x,y,\lambda \right)=0$ or $\|\nabla G^{\alpha,\beta,\gamma}\left( x,y,\lambda \right)\|\le \varepsilon$,  we call the point $(x,y)$ is a stationary point or an $\varepsilon$-stationary point of problem \eqref{problem:1}, respectively.
\end{definition}

Note that a stationary point of problem \eqref{problem:1} is also a KKT point of problem \eqref{problem:1}. Moreover, $\|\nabla G^{\alpha,\beta,\gamma}(x,y,\lambda)\|\leq\varepsilon$ implies that  $\|Ax+By-c\|\leq \varepsilon$, which means that an $\varepsilon$-stationary point of problem \eqref{problem:1} is also an $\varepsilon$-KKT point of problem \eqref{problem:1}. We refer to \cite{Dussault,Kanzow,lu23} for more details of the KKT point and  the $\varepsilon$-KKT point.

\subsection{A Primal-Dual Alternating Proximal Gradient Algorithm}
	
By the strong duality shown in Theorem \ref{dual}, instead of solving \eqref{problem:1}, we propose a primal-dual alternating proximal gradient (PDAPG) algorithm for solving \eqref{problem:2}. 
At each iteration of the proposed PDAPG algorithm, it performs a proximal gradient step for updating $y$ and $x$, respectively and a gradient projection step for updating $\lambda$ for a regularized version of $\mathcal{L}(x,y,\lambda)$, i.e., 
\begin{align}
	\mathcal{L}_k(x,y,\lambda)&=\mathcal{L}(x,y,\lambda)-\frac{\rho_k}{2}\|y\|^2,\label{sc:1}
\end{align}
where $\rho_k\ge0$ is a regularization parameter.
More detailedly, at the $k$th iteration of the proposed algorithm, the update for $y_k$ is as follows,
\begin{align}
	y_{k+1}&={\arg\max}_{y\in \mathcal{Y}}\langle  \nabla _y \mathcal{L}_k\left( x_{k},y_k,\lambda_k \right),y-y_k \rangle  -\frac{\beta}{2}\| y-y_k \| ^2-g(y),
\end{align}
where $\beta >0$ is a parameter which will be defined later. The update for $x_k$ is as follows,
\begin{align}
	x_{k+1}&={\arg\min}_{x\in \mathcal{X}}\langle  \nabla _x \mathcal{L}_k\left( x_{k},y_{k+1},\lambda_k \right) ,x-x_k \rangle  +\frac{\alpha_k}{2}\| x-x_k \| ^2+h(x),
\end{align}
where $\alpha_k > 0$ is a parameter which will be defined later. The update for $\lambda_k$ is as follows,
\begin{align}
	\lambda_{k+1}&=\lambda_k-\gamma_k\nabla_{\lambda}\mathcal{L}(x_{k+1},y_{k+1},\lambda_k).
\end{align}
The proposed PDAPG algorithm is formally stated in Algorithm \ref{alg:1}.
\begin{algorithm}[t]
	\caption{(PDAPG)}
	\label{alg:1}
	\begin{algorithmic}
		\STATE \textbf{Step 1}:Input $x_1,y_1,\lambda_1,\beta,\alpha_1,\gamma_1$; Set $k=1$.
		\STATE \textbf{Step 2}:Perform the following update for $y_k$:  	
		\begin{equation}\label{update-y}
			y_{k+1}=\operatorname{Prox}_{g,\mathcal{Y}}^{\beta}\left( y_k+\frac{1}{\beta} \nabla _y \mathcal{L}\left( x_{k},y_k,\lambda_k \right) -\frac{\rho_k}{\beta} y_k\right).
		\end{equation}	
		\STATE \textbf{Step 3}:Perform the following update for $x_k$:
		\begin{equation}\label{update-x}
			x_{k+1}=\operatorname{Prox}_{h,\mathcal{X}}^{\alpha_k} \left( x_k - \frac{1}{{\alpha_k}}  \nabla _x\mathcal{L}\left( x_{k},y_{k+1},\lambda_k \right) \right).
		\end{equation}
		\STATE \textbf{Step 4}:Perform the following update for $\lambda_k$:
		\begin{equation}\label{update-lambda}
			\lambda_{k+1}=\lambda_k+\gamma_k(Ax_{k+1}+By_{k+1}-c).
		\end{equation}
		\STATE \textbf{Step 5}:If $\|\nabla G^{\alpha_k,\beta,\gamma_k}\left( x_{k+1},y_{k+1},\lambda_{k+1} \right)\|\le\varepsilon$, stop; otherwise, set $k=k+1, $ go to Step 2.
	\end{algorithmic}
\end{algorithm}

Note that PDAPG is a single-loop algorithm. When $f(x,y)$ is a general concave function with respect to $y$, PDAPG is a brand-new algorithm. We set $\rho_k=0$ when $f(x,y)$ is strongly concave with respect to $y$. Even for nonconvex-strongly concave minimax problem, PDAPG is different from PGmsAD in \cite{Dai} and MGD in  \cite{Tsaknakis} that are both multi-loop algorithms. Moreover, $x_{k+1}$ is used to update $\lambda_k$ in Algorithm \ref{alg:1} instead of $x_k$ which is used in PGmsAD algorithm.  

Note that for algorithmic aspects, under the general nonconvex-concave setting, when updating $y$, PDAPG updates $y$ via gradient projection or proximal gradient steps, while HiBSA in \cite{Lu} needs to maximize the original function $f (x, y)$ plus some regularization terms, this subproblem may be more difficult to solve in many cases. In addition, although problem \eqref{problem:2} is a non-smooth nonconvex-concave minimax problem, PDAPG cannot be regarded as a special case of BAPG \cite{Xu} since that $\Lambda$ in \eqref{problem:2} is not a compact set, which is a necessary condition to ensure the convergence of BAPG. This is also the main difficulty in analyzing the PDAPG algorithm. In order to overcome this difficulty, we need to construct a new potential function.

To analyze the iteration complexity of Algorithm \ref{alg:1}, we need to make the following assumption about the smoothness of $f(x,y)$.
\begin{assumption}\label{ass:Lip}
	The function $f(x,y)$ has Lipschitz continuous gradients and there exists a constant $L>0$ such that for any $x, x_1, x_2\in \mathcal{X}$, and $y, y_1, y_2\in \mathcal{Y}$, we have
	\begin{align*}
		\| \nabla_{x} f\left(x_1, y\right)-\nabla_{x} f\left(x_2, y\right)\| &\leq L\|x_{1}-x_{2}\|,\\
		\|\nabla_{x} f\left(x, y_1\right)-\nabla_{x} f\left(x, y_2\right)\| &\leq L\|y_{1}-y_{2}\|,\\
		\|\nabla_{y} f\left(x, y_1\right)-\nabla_{y} f\left(x, y_2\right)\| &\leq L\|y_{1}-y_{2}\|,\\
		\|\nabla_{y} f\left(x_1, y\right)-\nabla_{y} f\left(x_2, y\right)\| &\leq L\|x_{1}-x_{2}\|.
	\end{align*}
\end{assumption}

\subsection{Nonconvex-Strongly Concave Setting}
In this subsection, we prove the iteration complexity of Algorithm \ref{alg:1} under the nonconvex-strongly concave setting, i.e., $f(x,y)$ is $\mu$-strongly concave with respect to $y$ for any given $x\in\mathcal{X}$. Under this setting, $\forall k \ge 1$, we set
$
\alpha_k=\alpha,\gamma_k=\gamma, \rho_k=0.
$
Denote   
$\Phi(x,\lambda):=\max_{y\in\mathcal{Y}}\mathcal{L}(x,y,\lambda)$ and 
$y^*(x,\lambda):=\arg\max_{y\in\mathcal{Y}}\mathcal{L}(x,y,\lambda)$.
$\Phi(x,\lambda)$ is $L_\Phi$-Lipschitz smooth with $L_\Phi=L+\frac{L^2}{\mu}$ under the $\mu$-strong concavity of $\mathcal{L}(x,y,\lambda)$ by \cite{Nouiehed}. Moreover,  by $\mu$-strong concavity of $\mathcal{L}(x,y,\lambda)$ with respect to $y$,  for any given $x$ and $\lambda$ we have
\begin{align}
	\nabla_x\Phi(x,\lambda)=\nabla_x \mathcal{L}(x, y^*(x,\lambda),\lambda),\quad
	\nabla_\lambda\Phi(x,\lambda)=\nabla_\lambda \mathcal{L}(x, y^*(x,\lambda),\lambda)\label{gradla:nsc}.
\end{align}
\begin{lemma}\label{lem:eb}
	Suppose that Assumption \ref{ass:Lip} holds and $f(x,y)$ is $\mu$-strongly concave with respect to $y$. Let $\eta=\frac{(2\beta+\mu)(\beta+L)}{\mu\beta}$. Then
	\begin{align}\label{lem:eb:1}
		\|y-y^*(x,\lambda)\|\le\eta\left\| y-\operatorname{Prox}_{g,\mathcal{Y}}^{\beta}\left( y+ \frac{1}{\beta} \nabla _y\mathcal{L}\left( x,y,\lambda \right) \right)\right\|.
	\end{align}
\end{lemma}
\begin{proof}{Proof}
	Firstly, $g(y)$ is a closed convex function and $-\mathcal{L}(x,y,\lambda)$ is a convex function with $L$-Lipschitz continuous gradient with respect to $y$. Secondly, by the strongly convexity of $\varphi(x, y, \lambda):=-\mathcal{L}(x,y,\lambda)+g(y)$, it has a nonempty set $S$ of minimizers with respect to $y$ in $\mathcal{Y}$ and satisfies the quadratic growth condition $\varphi(x, y, \lambda)\ge\varphi^*+\frac{\mu}{2}dist^2(y,S)$, where $dist(y,S):=\inf_{z\in S}\|z-y\|$ and $\varphi^*=\min_{y\in \mathcal{Y}} \varphi(x, y, \lambda)$. Thus, by Corollary 3.6 in \cite{Drusvyatskiy}, the proof is completed. 
\end{proof}

Next, we provide an upper bound estimate of the difference between $\Phi(x_{k+1},\lambda_{k+1})$ and $\Phi(x_k,\lambda_k)$. All the proofs of the following lemmas and theorems in this section are shown in Appendix \ref{ap:sc} and Appendix \ref{ap:c}  for clearer reading.
\begin{lemma}\label{lem1}
Suppose that Assumption \ref{ass:Lip} holds. Let $\{\left(x_k,y_k,\lambda_k\right)\}$ be a sequence generated by Algorithm \ref{alg:1}. Then $\forall k \ge 1$, and for any $a>0$, $b>0$,
\begin{align}\label{lem1:iq1}
&\Phi(x_{k+1},\lambda_{k+1}) -\Phi(x_k,\lambda_k)\nonumber \\
\leq& \langle  \nabla _x\mathcal{L}( x_k,y_{k+1},\lambda_k),x_{k+1}-x_k \rangle+\frac{(L+\beta)^2\eta^2}{2\beta^2}\left(\frac{L^2}{a}+\frac{2\|B\|^2}{b} \right)\|y_{k+1}-y_k\|^2 \nonumber \\
&+\frac{ b+L_\Phi}{2}\|\lambda_{k+1}-\lambda_k\|^2+\left(\frac{\|B\|^2L^2\eta^2}{b\beta^2}+\frac{a+L_\Phi}{2}\right)\|x_{k+1}-x_k\|^2\nonumber\\
&+\langle \nabla _\lambda \mathcal{L}(x_{k+1},y_{k+1},\lambda_k) ,\lambda_{k+1}-\lambda_k \rangle.
\end{align}
\end{lemma}


Next, we provide an lower bound for the difference between $\mathcal{L}(x_{k+1},y_{k+1},\lambda_{k+1})$ and $\mathcal{L}(x_k,y_k,\lambda_k)$.
\begin{lemma}\label{lem2}
Suppose that Assumption \ref{ass:Lip} holds. Let $\{\left(x_k,y_k,\lambda_k\right)\}$ be a sequence generated by Algorithm \ref{alg:1}.
Then $\forall k \ge 1$,
\begin{align}\label{lem2:iq1}
	&\mathcal{L}(x_{k+1},y_{k+1},\lambda_{k+1})-\mathcal{L}( x_{k},y_{k},\lambda_k)\nonumber \\
	\geq&\left(\beta-\frac{L}{2}\right)\|y_{k+1}-y_k\|^2-\frac{L}{2}\| x_{k+1}-x_{k} \|^2+\langle \nabla _\lambda \mathcal{L}(x_{k+1},y_{k+1},\lambda_k),\lambda_{k+1}-\lambda_{k} \rangle\nonumber\\
	&+\langle\xi_{k+1},y_{k+1}-y_k \rangle+\langle \nabla _{x}\mathcal{L}(x_{k},y_{k+1},\lambda_k),x_{k+1}-x_{k} \rangle,
\end{align}
where $\xi_{k+1}\in\partial g(y_{k+1})$.

\end{lemma}

We now establish an important recursion for the Algorithm \ref{alg:1}.
\begin{lemma}\label{lem3}
Suppose that Assumption \ref{ass:Lip} holds. Denote
$
S(x,y,\lambda)=2\Phi(x,\lambda)-\mathcal{L}(x,y,\lambda)+h(x)+g(y).
$
Let $\{\left(x_k,y_k,\lambda_k\right)\}$ be a sequence generated by Algorithm \ref{alg:1}. Then $\forall k \geq 1$,
\begin{align}\label{lem3:1}
&S(x_{k+1},y_{k+1},\lambda_{k+1})-S( x_{k},y_{k},\lambda_k)\nonumber \\
\le&-\left(\alpha-\frac{L^3}{(L+\beta)^2}-\frac{L(L+\beta)^2\eta^2}{\beta^2}
-\frac{L^2}{\mu}-\frac{3L}{2}\right)\|x_{k+1}-x_k\|^2 \nonumber \\
&-\left(\frac{1}{\gamma}-\frac{2\|B\|^2(L+\beta)^2\eta^2}{L\beta^2}-L-\frac{L^2}{\mu} \right)\|\lambda_{k+1}-\lambda_{k}\|^2-\left(\beta-\frac{5L}{2}\right)\|y_{k+1}-y_k\|^2.
\end{align}
\end{lemma}

Let $\nabla G^{\alpha,\beta,\gamma}( x_k,y_k,\lambda_k)$ be defined as Definition \ref{gap-f}, we provide an upper bound on $\| \nabla G^{\alpha,\beta,\gamma}(x_k,y_k,\lambda_k)\|$ in the following lemma.

\begin{lemma}\label{lem4}
Suppose that Assumption \ref{ass:Lip} holds. Let $\{\left(x_k,y_k,\lambda_k\right)\}$ be a sequence generated by Algorithm \ref{alg:1}. Then $\forall k \geq 1$,
\begin{align}\label{lem4:1}
&\|\nabla G^{\alpha,\beta,\gamma}( x_k,y_{k},\lambda_k)\|^2
\le (\beta^2+2L^{2}+3\|B\|^2)\|y_{k+1}-y_k\|^2+(2\alpha^2+3\|A\|^2)\|x_{k+1}-x_k\|^2\nonumber\\
&+\frac{3}{\gamma^2}\|\lambda_{k+1}-\lambda_k\|^2.
\end{align}
\end{lemma}

Denote $T(\varepsilon):=\min\{k \mid \|\nabla G^{\alpha,\beta,\gamma}(x_k,y_k,\lambda_k)\|\leq \varepsilon \}$ and $\varepsilon>0$ is a given target accuracy. The following theorem gives a bound on $T(\varepsilon)$. 

\begin{theorem}\label{thm1}
Suppose that Assumption \ref{ass:Lip} hold. Let $\{\left(x_k,y_k,\lambda_k\right)\}$ be a sequence generated by Algorithm \ref{alg:1}. Let $\eta=\frac{(2\beta+\mu)(\beta+L)}{\mu\beta}$. If
\begin{align}\label{thm1:1}
&\beta>\frac{5L}{2}, \alpha>\frac{L^3}{(L+\beta)^2}+\frac{L(L+\beta)^2\eta^2}{\beta^2}
+\frac{L^2}{\mu}+\frac{3L}{2},\frac{1}{\gamma}>\left(\frac{2\|B\|^2(L+\beta)^2\eta^2}{L\beta^2}+L+\frac{L^2}{\mu} \right),
\end{align}
then $\forall \varepsilon>0$, we have
$$ T\left( \varepsilon \right) \le \frac{d_2}{\varepsilon ^2 d_1},$$
where $d_1:=\frac{\min\{\alpha-\frac{L^3}{(L+\beta)^2}-\frac{L(L+\beta)^2\eta^2}{\beta^2}
-\frac{L^2}{\mu}-\frac{3L}{2},\ \beta-\frac{5L}{2},\ \frac{1}{\gamma}-\frac{2\|B\|^2(L+\beta)^2\eta^2}{L\beta^2}-L-\frac{L^2}{\mu}\} }{\max \{ \beta^2+2L^{2}+3\|B\|^2,\ 2\alpha^2+3\|A\|^2,\ 3/\gamma^2\} }$, $d_2:=S(x_{1},y_{1},\lambda_{1})-\underbar{S}$ with $\underbar{S}:=\min\limits_{\lambda,x\in\mathcal{X},y\in\mathcal{Y}}S(x,y,\lambda)$.
\end{theorem}
%

\begin{remark}
Denote $\kappa=L/\mu$. If we set $\beta=3L$, by Theorem \ref{thm1} we have that $\eta=\mathcal{O}(\kappa)$, $\alpha=\mathcal{O}(\kappa^2)$ and $1/\gamma=\mathcal{O}(\kappa^2)$, which implies that the number of iterations for Algorithm \ref{alg:1} to obtain an $\varepsilon$-stationary point of \eqref{problem:1} is bounded by $\mathcal{O}\left(\kappa^2\varepsilon ^{-2} \right)$ under the nonconvex-strongly concave setting.
\end{remark}

\begin{remark}
In the strongly concave case, we can replace the assumption of compacity on $\mathcal{X}$ and $\mathcal{Y}$ by the more general assumption that the primal function $\max\limits_{\substack{y\in\mathcal{Y}\\Ax+By= c}}\{f(x,y)+h(x)-g(y)\} $ is lower-bounded.
\end{remark}

\subsection{Nonconvex-Concave Setting}
In this subsection, we prove the iteration complexity of Algorithm \ref{alg:1} under the nonconvex-concave setting. $\forall k\ge 1$, we first denote
\begin{align}
\Psi_k(x,\lambda)=\max_{y\in\mathcal{Y}}\mathcal{L}_k(x,y,\lambda),
\tilde{y}_k^*(x,\lambda)=\arg\max_{y\in\mathcal{Y}}\mathcal{L}_k(x,y,\lambda).\label{sc:3}
\end{align}
We also need to make the following assumption on the parameter $\rho_k$.
\begin{assumption}\label{muk}
$\{\rho_k\}$ is a nonnegative monotonically decreasing sequence.
\end{assumption}

\begin{lemma}\label{sc:lem1}
Suppose that Assumption \ref{ass:Lip} holds. Denote
$
M_k(x,y,\lambda)=2\Psi_k(x,\lambda)-\mathcal{L}_k(x,y,\lambda)+h(x)+g(y).
$
Let $\{\left(x_k,y_k,\lambda_k\right)\}$ be a sequence generated by Algorithm \ref{alg:1}. Setting $\eta_k=\frac{(2\beta+\rho_k)(\beta+L)}{\rho_k\beta}$, then $\forall k \geq 1$,
\begin{align}
&M_{k+1}(x_{k+1},y_{k+1},\lambda_{k+1})-M_k( x_{k},y_{k},\lambda_k)\nonumber \\
\le&-\left(\alpha_k-\frac{L^3}{(L+\beta)^2}-\frac{L(L+\beta)^2\eta_k^2}{\beta^2}
-\frac{L^2}{\rho_k}-\frac{3L}{2}\right)\|x_{k+1}-x_k\|^2 -\left(\beta-\frac{5L}{2}\right)\|y_{k+1}-y_k\|^2\nonumber \\
&-\left(\frac{1}{\gamma_k}-\frac{2\|B\|^2(L+\beta)^2\eta_k^2}{L\beta^2}-L-\frac{L^2}{\rho_k} \right)\|\lambda_{k+1}-\lambda_{k}\|^2+(\rho_k-\rho_{k+1})\sigma_y^2,\label{sc:lem1:1}
\end{align}
where $\sigma_y=\max\{\|y\| \mid y\in\mathcal{Y}\}$.
\end{lemma}

Similar to the proof of Theorem \ref{thm1}, we then obtain the following theorem which provides a bound on $T(\varepsilon)$, where $T(\varepsilon):=\min\{k \mid \|\nabla G^{\alpha_k,\beta,\gamma_k}(x_k,y_k,\lambda_k)\|\leq \varepsilon \}$ and $\varepsilon>0$ is a given target accuracy.

\begin{theorem}\label{sc:thm1}
Suppose that Assumptions \ref{ass:Lip} and \ref{muk} hold. Let $\{\left(x_k,y_k,\lambda_k\right)\}$ be a sequence generated by Algorithm \ref{alg:1}. Set  $\alpha{\tiny }_k=\frac{L^3}{(L+\beta)^2}+\frac{L\tau(L+\beta)^4(2\beta+\rho_k)^2}{\beta^4\rho_k^2}
+\frac{\tau L^2}{\rho_k}+\frac{3L}{2}$, $\frac{1}{\gamma_k}=\frac{[2\|B\|^2+L^2(\tau-1)](L+\beta)^4(2\beta+\rho_k)^2}{L\beta^4\rho_k^2}+L+\frac{L^2\tau}{\rho_k}$ with $\rho_k=\frac{2(L+\beta)}{k^{1/4}}$, and $\beta>\frac{5L}{2}$, $\forall k\geq 1$. Then for any given $\varepsilon>0$,
$$T( \varepsilon)\le \max\left\{\frac{4\tilde{d}_3^2\tilde{d}_4^2L^2(\tau-1)^2[2(L+\beta)^3(L+2\beta)^2+L\beta^4]^2}{\beta^8(L+\beta)^2\varepsilon^4},\frac{256(L+\beta)^4\sigma_y^4}{\varepsilon^4}\right\},$$
where $\tilde{d}_1=\frac{8\tau^2}{(\tau-1)^2} + \frac{\frac{8L^6\beta^4}{(L+\beta)^4}+12L^2\beta^4+3\|A\|^2\beta^4}{L^2(L+\beta)^4(\tau-1)^2}+\frac{2\tau^2L^4\beta^4}{L^2(L+\beta)^6(\tau-1)^2}$, $\tilde{d}_2=\frac{9[2\|B\|^2+L^2(\tau-1)]^2}{L^2(\tau-1)^4} + \frac{9L^2\beta^4((L+\beta)^2+\tau^2L^2)}{L^2(L+\beta)^6(\tau-1)^2}$,  $\tilde{d}_3=M_{1}(x_{1},y_{1},\lambda_1)-\underline{M}+\rho_1\sigma_y^2$ with $\underline{M}:=\min\limits_{\lambda\in\Lambda}\min\limits_{x\in\mathcal{X}}\min\limits_{y\in\mathcal{Y}}M_k(x,y,\lambda)$, $\tilde{d}_4=\max\left\{\frac{(\beta^2+2L^{2}+3\|B\|^2)\rho_1}{(\beta-\frac{5L}{2})L^2(\tau-1)}, \tilde{d}_1, \tilde{d}_2\right\}$ and $\tau>1$.
\end{theorem}

\begin{remark}
Theorem \ref{sc:thm1} implies that the number of iterations for Algorithm \ref{alg:1} to obtain an $\varepsilon$-stationary point of \eqref{problem:1} is bounded by $\mathcal{O}\left(\varepsilon ^{-4} \right)$, which matches the best known iteration complexity of single loop algorithms under the nonconvex-concave setting. The optimal iteration complexity of the existing multi-loop algorithms for solving smooth nonconvex-concave minimax problems is $\tilde{\mathcal{O}}\left( \varepsilon ^{-2.5} \right)$. Although we cannot theoretically prove that the iteration complexity of Algorithm \ref{alg:1} can reach this convergence order currently, we find that the numerical performance of Algorithm \ref{alg:1} is consistent with that of the two existing multi-loop algorithms, MGD and PGmsAD in Section 4.
\end{remark}

\section{Numerical Results}
In this section, we test the proposed PDAPG algorithm with MGD \cite{Tsaknakis} algorithm, PGmsAD \cite{Dai} algorithm for solving four test problems: an adversarial attack in network flow problem, a generalized absolute value equation, a linear regression problem, and a distributionally robust optimization problem.
All these problems are formulated as minimax problems with coupled constraints.

\subsection{Adversarial attack in network flow problems}
We consider the adversarial attack in network flow problem \eqref{testprobn1}. We generate networks (graphs) of $n$ nodes at random using the Erdos-Renyi model with parameter $p$ (i.e., the probability that an edge appears on the graph is $p$).

\begin{figure}[htb]
\centering  
\subfigure[n=15, p=0.75, d=10]{
	\includegraphics[width=0.4\textwidth]{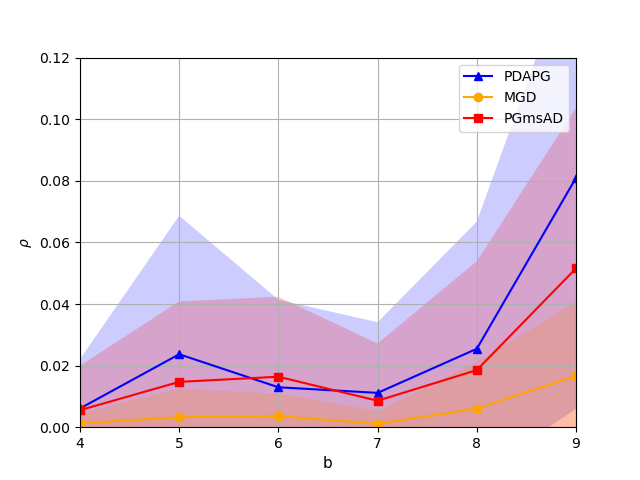}}
\subfigure[n=15, p=1, d=10]{
	\includegraphics[width=0.4\textwidth]{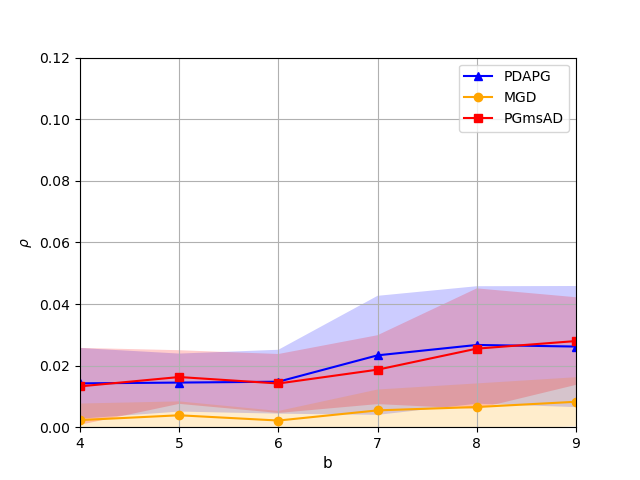}}
\caption{Relative total cost increase of three algorithms with $d=10$.}
\label{Fig2}
\end{figure}
\begin{figure}[htb]
\centering 
\subfigure[n=15, p=0.75, d=20]{
	\includegraphics[width=0.4\textwidth]{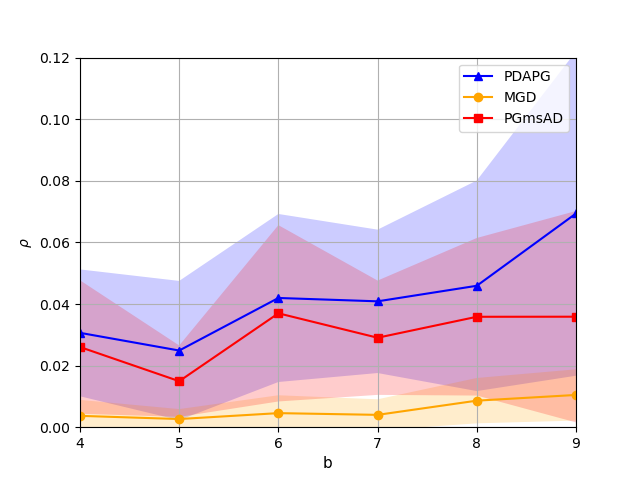}}
\subfigure[n=15, p=1, d=20]{
	\includegraphics[width=0.4\textwidth]{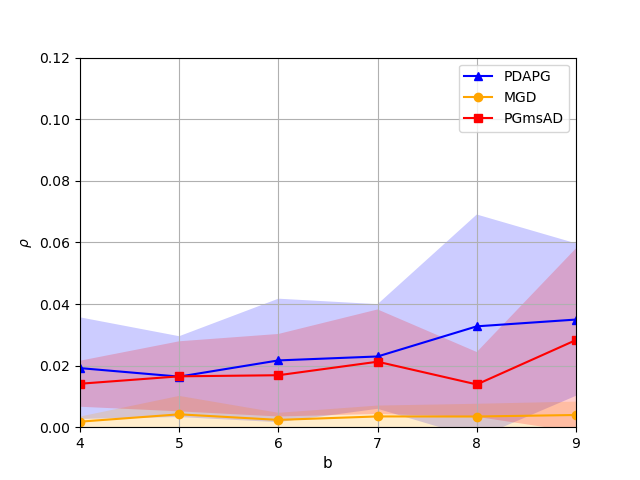}}
\caption{Relative total cost increase of three algorithms with $d=20$.}
\label{Fig3}
\end{figure}

We set $\eta=0.05$,  $1/\beta=0.8$, $1/\alpha=0.6$, $\gamma=0.5$ in PDAPG. All the stepsizes for MGD are set to be 0.5 which is the same as that in \cite{Tsaknakis}. We set $\alpha_y=0.8$, $\alpha_x=0.6 $ for PGmsAD algorithm which is the same as that of PDAPG algorithm. Following the same setting as that in \cite{Tsaknakis}, the capacity $p_{ij}$ and the cost coefficients $q_{ij}$ of the edges are generated uniformly at random in the interval [1, 2],  and total flow $r_t$ is $d\%$ of the sum of capacities of the edges exiting the source with $d$ being a parameter to be chosen. 

Following the same performance measure as that in \cite{Tsaknakis}, we also use relative increase of the cost $\rho=\frac{q_{att}-q_{cl}}{q_{cl}}$ as the performance measure, 
where $q_{cl}$ and $q_{att}$ denote the minimum cost before and after the attack, respectively. The higher the increase of the minimum cost, the more successful the attack. We compare the relative total cost increase of the three tested algorithms for solving \eqref{testprobn1} with different $b$. For each given $b$, all three algorithms run 15 times with randomly chosen $p_{ij}$ and $q_{ij}$ which means 15 graphs with different capacity and cost vectors.

Fig. \ref{Fig2}-\ref{Fig3} show the average of the relative total cost increase of three tested algorithms over 15 independent runs with different $d$. The shaded part around lines denotes the standard deviation over 15 independent runs.
We can find that the proposed PDAPG algorithm slightly outperforms the MGD algorithm, and is similar to that of the PGmsAD algorithm.

\subsection{Generalized absolute value equations}
In this subsection, we consider the generalized absolute value equations  \eqref{testprob2} discussed in Section \ref{app}.

We compared PDAPG with PGmsAD and MGD for solving generalized absolute value equations. Following the same setting as that in \cite{Dai}, we set
\begin{equation}
	A=\left(\begin{array}{lll}
		1 & 1 & 1 \\
		1 & 0 & 1 \\
		1 & 1 & 1
	\end{array}\right), \quad B=\left(\begin{array}{ccc}
		-1 & 1& 0 \\
		1 & 2 & 1 \\
		0 & 1 & 1
	\end{array}\right), \quad b=(-1,4,1)^{\top} .
\end{equation}
Initial points are chosen as zero vectors in all algorithms. In both PGmsAD and MGD, the number of inner loops was selected as $N = 5$. For the PDAPG algorithm, we set $\rho_k=k^{-1/4}$, $1/\beta=0.002$, $1/\alpha=0.01k^{-1/2}$, $\gamma=0.01k^{-1/2}$.
All the stepsizes for MGD are set to be 0.001. We set $\alpha_y=0.003$, $\alpha_x=0.002$ for PGmsAD algorithm. We use error 
$\|Ax + B|x| - b\|$ defined in \eqref{gave} as the performance measure, which is same as that in \cite{Dai}.

\begin{figure}[htb]
	\centering 
	\includegraphics[width=0.5\textwidth]{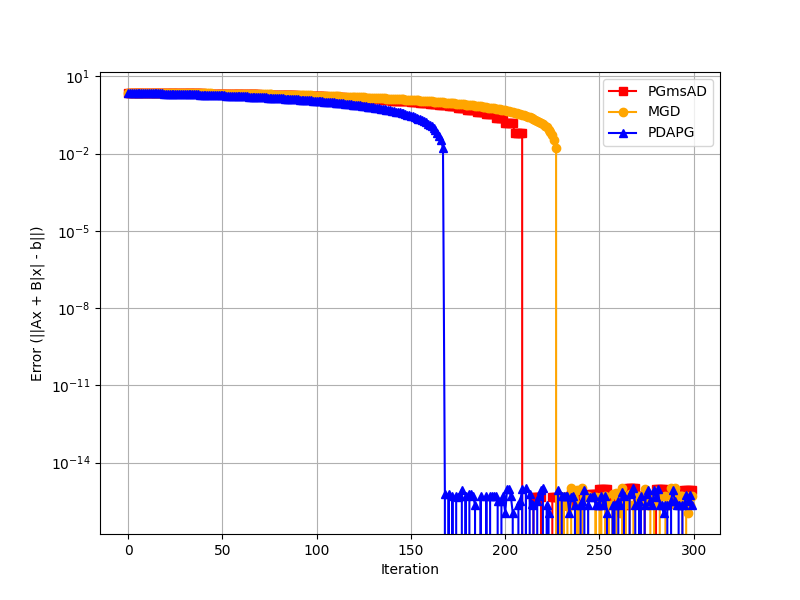}
	\caption{Error curves for generalized absolute value equation of three algorithms.}
	\label{Fig4}
\end{figure}

Fig. \ref{Fig4} shows the error curves for generalized absolute value equation of three tested algorithms.
We can find that the proposed PDAPG algorithm slightly outperforms the MGD algorithm and PGmsAD algorithm.

\subsection{Linear regression:}
In this subsection, we consider the well-known linear regression problem with joint linear constraints as follows \cite{Dai}:
\begin{align}\label{test-lr}
	\min _{x \in \mathbb{R}^n} \max _{y \in \mathbb{R}^m}  &f(x,y)=\frac{1}{m}[-\frac{1}{2}\|y\|^2-b^{\top}y+y^{\top}Kx]+\frac{1}{2m}\|x\|^2\nonumber\\
	\text {s.t.} \quad&Ax+By+c=0,
\end{align}
where the matrix $K \in \mathbb{R}^{m\times n}$, $A \in  \mathbb{R}^{p\times n}$, and $B \in  \mathbb{R}^{p\times m}$ are generated by a Gaussian distribution $ N (0, I)$.

We compared PDAPG with PGmsAD and MGD for solving linear regression. Let $n =m$, $b=0$, $c=0$, the initial points were randomly selected in all algorithms. In both PGmsAD and MGD, the number of inner loops was selected as $N = 3$, which is same as that in \cite{Dai}.
We set the step sizes $1/\beta=0.1$, $1/\alpha=0.001$, $\gamma=0.001$ in PDAPG. For the MGD algorithm, we set $\beta=0.1$, $\alpha=0.001$. For the PGmsAD algorithm, we set $\alpha_y=0.1$, $\alpha_x=0.001$.

\begin{figure}[htb]
	\centering  
	\subfigure[n=100, p=10]{
		\includegraphics[width=0.45\textwidth]{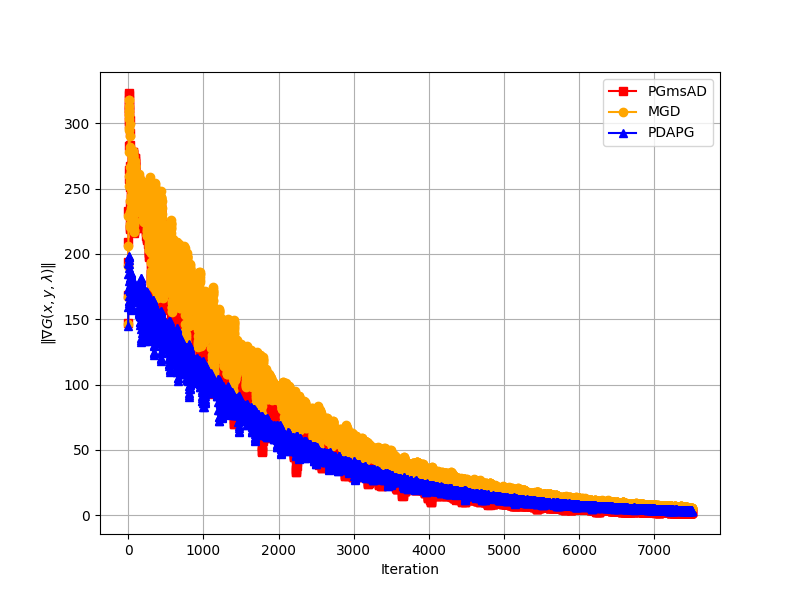}}
	\subfigure[n=100, p=40]{
		\includegraphics[width=0.45\textwidth]{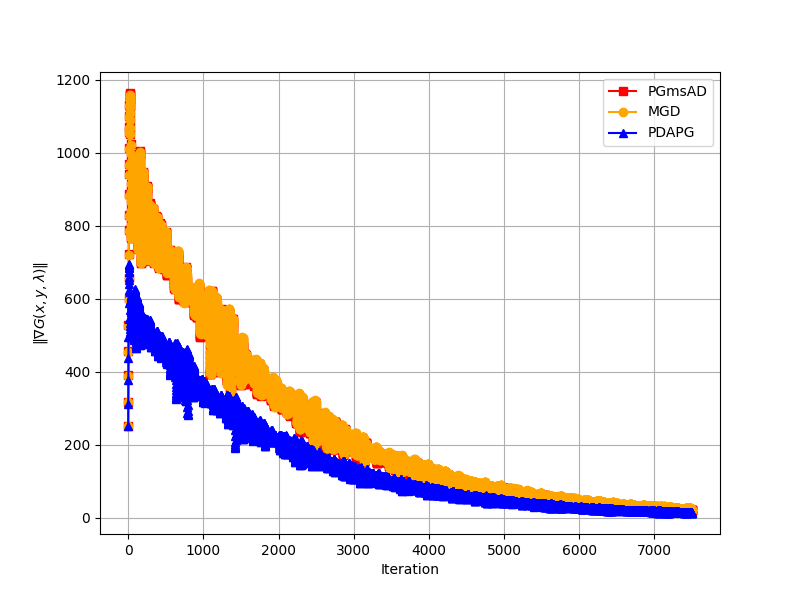}}
	\caption{Numerical results of the three tested algorithms for solving linear regression problem \eqref{test-lr}.}
	\label{Fig5}                                                         
\end{figure}

Fig. \ref{Fig5} shows the trend of the gradient norm of three tested algorithms. $n$ refers to the dimension of $x$ and $y$, $p$ is the 
first dimension of the matrix $A$ and $B$. We can find that the proposed PDAPG algorithm slightly outperforms the MGD algorithm, and is similar to that of the PGmsAD algorithm  in sub-figure (a). In sub-figure (b), the proposed PDAPG algorithm slightly outperforms the MGD algorithm and PGmsAD algorithm.

\subsection{Distributionally robust optimization problems}
	In this subsection, we investigate the distributionally robust optimization problem \eqref{test-dro} discussed in Section \ref{sec1}, which presents a nonconvex-concave minimax structure with non-separable coupled constraints.
	We conduct a comparative study of our proposed PDAPG against two state-of-the-art methods: PGmsAD and MGD. The logistic loss 
$\ell_i(x)=\log \left(1+\exp \left(-b_i \mathbf{a}_i^{\top} \mathbf{x}\right)\right)$ is adopted, where $\left\{\mathbf{a}_i, b_i\right\}_{i=1}^n$ denotes the dataset with features $\mathbf{a}_i \in \mathbb{R}^d$ and binary labels $b_i \in\{-1,1\}$. To evaluate scalability in high-dimensional settings, experiments are performed on the gisette dataset from LIBSVM ($n = 6000$, $d=5000$), a benchmark characterized by sparse large-scale features.
	
	Algorithmic configurations are set as follows: For PDAPG, we employ time-varying step sizes  $1/\beta=0.1$, $1/\alpha_k=\frac{0.5}{k^{0.5}+100}$, $\gamma_k=\frac{0.0001}{k^{0.5}+1000}$, $\rho_k=\frac{0.01}{k^{0.25}+200}$ to adaptively balance primal-dual updates. Comparative methods use fixed step sizes (MGD: $\beta=0.1$, $\alpha=0.001$; PGmsAD: $\alpha_y=0.1$, $\alpha_x=0.009$), limiting their capability to handle coupling constraints.
	\begin{figure}[htb]
		\centering 
		\includegraphics[width=0.5\textwidth]{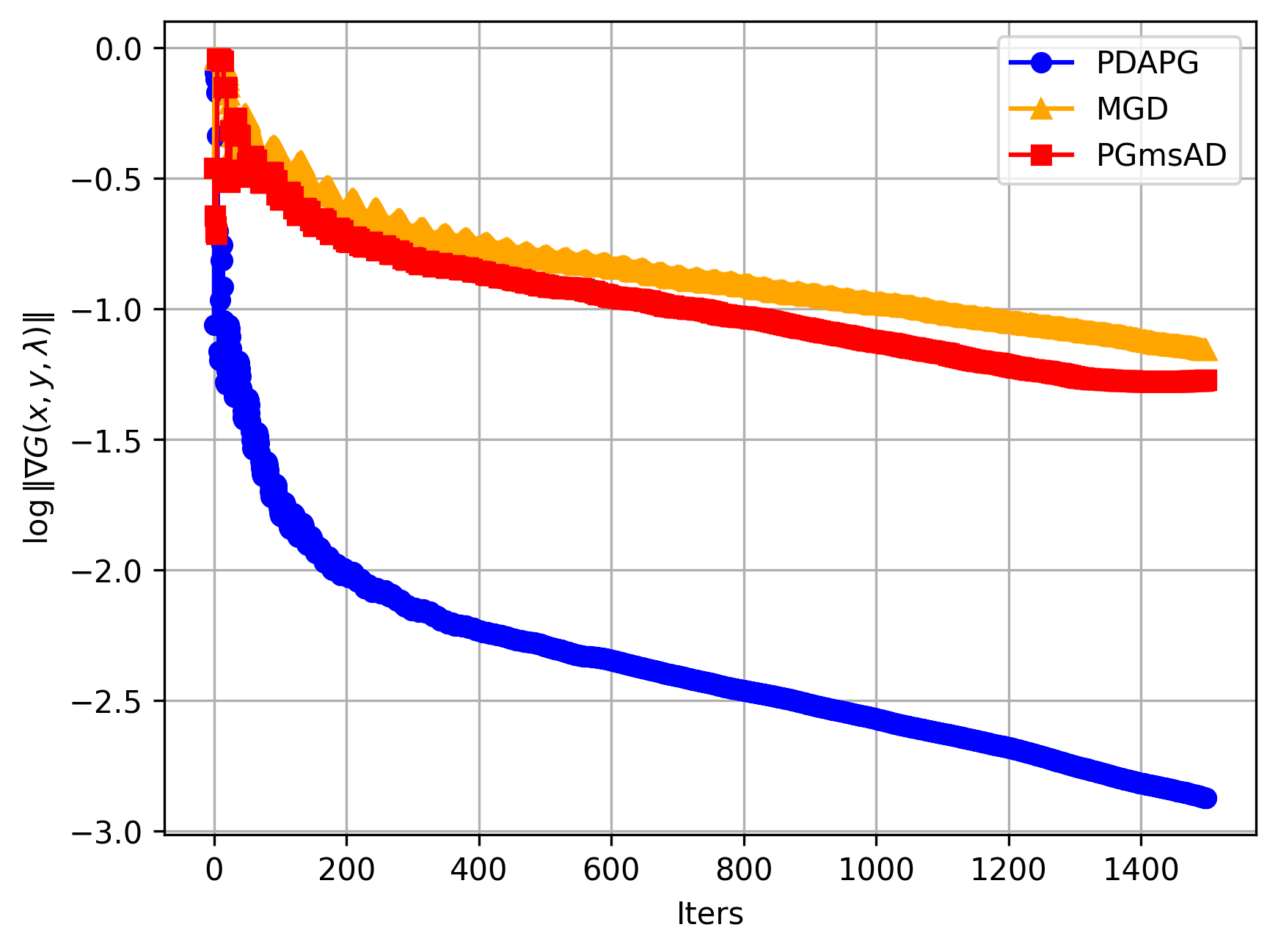}
		\caption{Numerical result of the PDAPG algorithm for solving distributionally robust optimization problem \eqref{test-dro}.}
		\label{Fig6}
	\end{figure}
	
	Convergence behavior is quantified through gradient norm trajectories in Fig. \ref{Fig6}. PDAPG achieves a decisive advantage, reducing the gradient norm much faster than PGmsAD and MGD. This empirically validates PDAPG's dual mechanism superiority, the decaying step-size scheme effectively coordinates primal-dual variables under coupling constraints.

\section{Conclusions}\label{section_conclu}
This work addresses the critical challenge of solving nonsmooth nonconvex-(strongly) concave minimax problems with coupled linear equality constraints—a class of iterative minimax games pervasive in adversarial learning, resource allocation, and strategic decision-making. Our contributions bridge theoretical and algorithmic gaps in constrained equilibrium-seeking dynamics. First, we establish the strong duality with respect to $y$ under feasible constraint qualifications, a foundational result enabling duality-based methods for nonconvex minimax games with linear coupling. Leveraging this duality, we propose the Primal-Dual Alternating Proximal Gradient (PDAPG) algorithm, a single-loop framework that iteratively updates primal ($x$) and dual $(y,\lambda)$ variables via proximal steps, effectively resolving strategic competition and constraint satisfaction in tandem.

Theoretical analysis demonstrates that PDAPG achieves $\mathcal{O}\left( \varepsilon ^{-2} \right)$  (nonconvex-strongly concave) and $\mathcal{O}\left( \varepsilon ^{-4} \right)$ (nonconvex-concave) iteration complexity to attain 
$\varepsilon$-stationary points—matching the lower-bound complexity for unconstrained nonconvex-strongly concave problems and aligning with state-of-the-art rates for unconstrained nonconvex-concave settings. Notably, PDAPG attains these guarantees while preserving a single-loop structure, circumventing nested iterations required by prior multi-loop methods like MGD \cite{Tsaknakis}. Numerical experiments validate the efficiency of our framework in practical scenarios.

While PDAPG advances the frontier of constrained minimax optimization, open questions remain. A key direction is whether improved complexity bounds can be achieved for nonsmooth nonconvex-(strongly) concave games with linear coupling, particularly through adaptive step sizes or problem structure exploitation. Extending this framework to stochastic or decentralized settings may further broaden its applicability in large-scale adversarial learning systems.

\section{Acknowledgments}  We thank Professor Guanghui (George) Lan of Georgia Institute of Technology for his suggestions that greatly improved the quality of this paper.

\section{Appendix} 
\subsection{The Proofs of Lemmas and Theorems in Section 2.3}\label{ap:sc}

\textbf{The proof of Lemma \ref{lem1}.}
\begin{proof}{Proof}
		\noindent
Since that $\Phi(x, \lambda)$ is $L_\Phi$-smooth with respect to $x$ and by \eqref{gradla:nsc}, we have that
\begin{align} \label{lem1:2}
	\Phi(x_{k+1},\lambda_{k})-\Phi(x_{k},\lambda_k)
				\le &\langle \nabla _x\Phi(x_{k},\lambda_k) ,x_{k+1}-x_k \rangle  +\frac{L_\Phi}{2}\|x_{k+1}-x_k\|^2\nonumber\\
	=&\langle \nabla _x\mathcal{L}(x_{k},y^*(x_k,\lambda_k),\lambda_k)-\nabla _x\mathcal{L}(x_{k},y_{k+1},\lambda_k) ,x_{k+1}-x_k \rangle  \nonumber\\
	&+\frac{L_\Phi}{2}\|x_{k+1}-x_k\|^2+\langle \nabla _x\mathcal{L}(x_{k},y_{k+1},\lambda_k) ,x_{k+1}-x_k \rangle.
\end{align}
By the Cauchy-Schwarz inequality, for any $a>0$, we have
\begin{align}
	&\langle \nabla _x\mathcal{L}(x_{k},y^*(x_k,\lambda_k),\lambda_k)-\nabla _x\mathcal{L}(x_{k},y_{k+1},\lambda_k) ,x_{k+1}-x_k \rangle\nonumber\\
	\le&\frac{L^2}{2a}\|y_{k+1}-y^*(x_k,\lambda_k)\|^2+\frac{a}{2}\|x_{k+1}-x_k\|^2.\label{lem1:3}
\end{align}
By \eqref{lem:eb:1} in Lemma \ref{lem:eb}, we further have
\begin{align}\label{lem1:4}
	&\|y_{k+1}-y^*(x_k,\lambda_k)\|
	\le \eta\left\| y_{k+1}-\operatorname{Prox}_{g,\mathcal{Y}}^{\beta}\left( y_{k+1}+ \frac{1}{\beta} \nabla _y\mathcal{L}\left( x_k,y_{k+1},\lambda_k \right) \right)\right\|.
\end{align}
By \eqref{update-y}, the non-expansive property of the proximity operator $\operatorname{Prox}_{g,\mathcal{Y}}^{\beta}(\cdot)$ and Assumption \ref{ass:Lip}, we obtain
\begin{align}\label{lem1:5}
	&\left\| y_{k+1}-\operatorname{Prox}_{g,\mathcal{Y}}^{\beta}\left( y_{k+1}+ \frac{1}{\beta} \nabla _y\mathcal{L}\left( x_k,y_{k+1},\lambda_k \right) \right)\right\|\nonumber\\
	= &\left\| \operatorname{Prox}_{g,\mathcal{Y}}^{\beta}\left( y_{k}+ \frac{1}{\beta} \nabla _yL\left( x_k,y_{k},\lambda_k \right) \right)-\operatorname{Prox}_{g,\mathcal{Y}}^{\beta}\left( y_{k+1}+ \frac{1}{\beta} \nabla _yL\left( x_k,y_{k+1},\lambda_k \right) \right)\right\|\nonumber\\
	\nonumber\\
	\le&\|y_{k+1}-y_k\|+\frac{1}{\beta}\|\nabla _yf\left( x_k,y_{k+1}\right)-\nabla _yf\left( x_k,y_{k} \right)\|
	\le\frac{L+\beta}{\beta}\|y_{k+1}-y_k\|.
\end{align}
Combing \eqref{lem1:2}, \eqref{lem1:3}, \eqref{lem1:4} and \eqref{lem1:5}, we get
\begin{align}\label{lem1:6}
	\Phi(x_{k+1},\lambda_{k})-\Phi(x_{k},\lambda_k)
	\le&\frac{L^2\eta^2(L+\beta)^2}{2a\beta^2}\|y_{k+1}-y_k\|^2+\frac{a+L_\Phi}{2}\|x_{k+1}-x_k\|^2\nonumber\\
&+\langle \nabla _x\mathcal{L}(x_{k},y_{k+1},\lambda_k) ,x_{k+1}-x_k \rangle.
\end{align}
On the other hand, $\Phi(x, \lambda)$ is $L_\Phi$-smooth with respect to $\lambda$ and \eqref{gradla:nsc}, we have
\begin{align} \label{lem1:7}
\Phi(x_{k+1},\lambda_{k+1})-\Phi(x_{k+1},\lambda_k)
\le&\langle \nabla _\lambda \mathcal{L}(x_{k+1},y^*(x_{k+1},\lambda_k),\lambda_k)-\nabla _\lambda \mathcal{L}(x_{k+1},y_{k+1},\lambda_k) ,\lambda_{k+1}-\lambda_k \rangle\nonumber\\
&+\langle \nabla _\lambda \mathcal{L}(x_{k+1},y_{k+1},\lambda_k) ,\lambda_{k+1}-\lambda_k \rangle+\frac{L_\Phi}{2}\|\lambda_{k+1}-\lambda_k\|^2.
\end{align}
Next, we estimate the first term in the right hand side of \eqref{lem1:7}. By the Cauchy-Schwarz inequality, for any $b>0$ we obtain
\begin{align} \label{lem1:8}
&\langle \nabla _\lambda \mathcal{L}(x_{k+1},y^*(x_{k+1},\lambda_k),\lambda_k)-\nabla _\lambda \mathcal{L}(x_{k+1},y_{k+1},\lambda_k) ,\lambda_{k+1}-\lambda_k \rangle\nonumber\\
\le&\frac{\|B\|^2}{2b}\|y_{k+1}-y^*(x_{k+1},\lambda_k)\|^2+\frac{b}{2}\|\lambda_{k+1}-\lambda_k\|^2.
\end{align}
By \eqref{lem:eb:1} in Lemma \ref{lem:eb}, we further have
\begin{align}
&\|y_{k+1}-y^*(x_{k+1},\lambda_k)\|
\le \eta\left\| y_{k+1}-\operatorname{Prox}_{g,\mathcal{Y}}^{\beta}\left( y_{k+1}+ \frac{1}{\beta} \nabla _y\mathcal{L}\left( x_{k+1},y_{k+1},\lambda_k \right) \right)\right\|.\label{lem1:9}
\end{align}
By \eqref{update-y}, the non-expansive property of the proximity operator $\operatorname{Prox}_{g,\mathcal{Y}}^{\beta}(\cdot)$ and Assumption \ref{ass:Lip}, we obtain
\begin{align}\label{lem1:10}
&\left\| y_{k+1}-\operatorname{Prox}_{g,\mathcal{Y}}^{\beta}\left( y_{k+1}+ \frac{1}{\beta} \nabla _y\mathcal{L}\left( x_{k+1},y_{k+1},\lambda_k \right) \right)\right\|
\le\frac{L+\beta}{\beta}\|y_{k+1}-y_k\|+\frac{L}{\beta}\|x_{k+1}-x_k\|.
\end{align}
By combing \eqref{lem1:7}, \eqref{lem1:8}, \eqref{lem1:9}, \eqref{lem1:10} and using the fact $(a+b)^2\le 2a^2+2b^2$, we get
\begin{align} \label{lem1:iq13}
\Phi(x_{k+1},\lambda_{k+1})-\Phi(x_{k+1},\lambda_k)
\le&\frac{\|B\|^2(L+\beta)^2\eta^2}{b\beta^2}\|y_{k+1}-y_k\|^2
+\frac{\|B\|^2L^2\eta^2}{b\beta^2}\|x_{k+1}-x_k\|^2\nonumber\\
&+\frac{b+L_\Phi}{2}\|\lambda_{k+1}-\lambda_k\|^2+\langle \nabla _\lambda \mathcal{L}(x_{k+1},y_{k+1},\lambda_k) ,\lambda_{k+1}-\lambda_k \rangle.
\end{align}
The proof is then completed by adding \eqref{lem1:6} and \eqref{lem1:iq13}.
%
\end{proof}

\textbf{The proof of Lemma \ref{lem2}.}
\begin{proof}{Proof}
		\noindent
The optimality condition for $y_k$ in \eqref{update-y} implies that $\forall y\in \mathcal{Y}$ and $\forall k\geq 1$,
\begin{equation}\label{lem2:3}
\langle \nabla _y\mathcal{L}(x_{k},y_k,\lambda_k),y_{k+1}-y_k \rangle \ge \beta\|y_{k+1}-y_k\|^2+\langle\xi_{k+1},y_{k+1}-y_k \rangle.
\end{equation}
By  Assumption \ref{ass:Lip}, the gradient of $\mathcal{L}(x,y,\lambda)$ is Lipschitz continuous with respect to $y$ and  \eqref{lem2:3}, we get
\begin{align}\label{lem2:4}
\mathcal{L}(x_{k},y_{k+1},\lambda_k)-\mathcal{L}( x_{k},y_{k},\lambda_k)
			\ge& \langle \nabla _{y}\mathcal{L}(x_{k},y_{k},\lambda_k),y_{k+1}-y_{k} \rangle -\frac{L}{2}\| y_{k+1}-y_{k} \|^2\nonumber\\
\ge&\left(\beta-\frac{L}{2}\right)\|y_{k+1}-y_k\|^2+\langle\xi_{k+1},y_{k+1}-y_k \rangle.
\end{align}
Since that the gradient of $\mathcal{L}(x,y,\lambda)$ is Lipschitz continuous with respect to $x$, we obtain
\begin{align}\label{lem2:5}
&\mathcal{L}(x_{k+1},y_{k+1},\lambda_k)-\mathcal{L}( x_{k},y_{k+1},\lambda_k)
\ge \langle \nabla _{x}\mathcal{L}(x_{k},y_{k+1},\lambda_k),x_{k+1}-x_{k} \rangle -\frac{L}{2}\| x_{k+1}-x_{k} \|^2.
\end{align}
On the other hand, it can be easily checked that
\begin{align}\label{lem2:iq6}
&\mathcal{L}(x_{k+1},y_{k+1},\lambda_{k+1})-\mathcal{L}( x_{k+1},y_{k+1},\lambda_k)
= \langle \nabla _\lambda \mathcal{L}(x_{k+1},y_{k+1},\lambda_k),\lambda_{k+1}-\lambda_{k} \rangle.
\end{align}
The proof is then completed by adding \eqref{lem2:4}, \eqref{lem2:5} and \eqref{lem2:iq6}.
\end{proof}

\textbf{The proof of Lemma \ref{lem3}.}
\begin{proof}{Proof}
		\noindent
The optimality condition for $\lambda_k$ in \eqref{update-lambda} implies that
\begin{equation}\label{lem3:iq4}
\langle \nabla _\lambda \mathcal{L}(x_{k+1},y_{k+1},\lambda_k)+\frac{1}{\gamma}(\lambda_{k+1}-\lambda_{k}),\lambda_k-\lambda_{k+1} \rangle \ge 0.
\end{equation}
Denote 	$S_1(x,y,\lambda)=2\Phi(x,\lambda)-\mathcal{L}(x,y,\lambda)$. 
Combining \eqref{lem1:iq1}, \eqref{lem2:iq1} and \eqref{lem3:iq4}, we have
\begin{align}\label{lem3:2}
&S_1(x_{k+1},y_{k+1},\lambda_{k+1})-S_1( x_{k},y_{k},\lambda_k)\nonumber \\
\le&-\left(\beta-\frac{L}{2}-\frac{(L+\beta)^2\eta^2}{\beta^2}\left(\frac{L^2}{a}+\frac{2\|B\|^2}{b} \right)\right)\|y_{k+1}-y_k\|^2 -\langle\xi_{k+1},y_{k+1}-y_k \rangle\nonumber \\
&+\left(\frac{2\|B\|^2L^2\eta^2}{b\beta^2}+a+L_\Phi+\frac{L}{2}\right)\|x_{k+1}-x_k\|^2+\langle  \nabla _x\mathcal{L}( x_k,y_{k+1},\lambda_k),x_{k+1}-x_k \rangle\nonumber\\
&-\left(\frac{1}{\gamma}- b-L_\Phi \right)\|\lambda_{k+1}-\lambda_k\|^2.
\end{align}
The optimality condition for $x_k$ in \eqref{update-x} implies that $\forall x\in \mathcal{X}$ and $\forall k\geq 1$,
\begin{equation}\label{lem3:4}
\langle \nabla _x\mathcal{L}(x_{k},y_{k+1},\lambda_k),x_{k+1}-x_k \rangle \le -\alpha\|x_{k+1}-x_k\|^2-\langle\zeta_{k+1},x_{k+1}-x_k \rangle.
\end{equation}
Plugging \eqref{lem3:4} into \eqref{lem3:2},
utilizing the convexity of $h(x)$ and $g(y)$, and by the definition of subgradient, we have $h(x_{k})-h( x_{k+1})\ge\langle\zeta_{k+1},x_{k}-x_{k+1} \rangle$ and $g(y_{k})-g( y_{k+1})\ge\langle\xi_{k+1},y_{k}-y_{k+1} \rangle$, 
then we obtain
\begin{align}\label{lem3:6}
&S(x_{k+1},y_{k+1},\lambda_{k+1})-S( x_{k},y_{k},\lambda_k)\nonumber \\
\le&-\left(\beta-\frac{L}{2}-\frac{(L+\beta)^2\eta^2}{\beta^2}\left(\frac{L^2}{a}+\frac{2\|B\|^2}{b} \right)\right)\|y_{k+1}-y_k\|^2 \nonumber \\
&-\left(\alpha-\frac{2\|B\|^2L^2\eta^2}{b\beta^2}-a-L_\Phi-\frac{L}{2}\right)\|x_{k+1}-x_k\|^2-\left(\frac{1}{\gamma}- b-L_\Phi \right)\|\lambda_{k+1}-\lambda_k\|^2.
\end{align}
The proof is then completed by $L_\Phi=L+\frac{L^2}{\mu}$ and choosing $a=\frac{L(L+\beta)^2\eta^2}{\beta^2}$ and $b=\frac{2\|B\|^2(L+\beta)^2\eta^2}{L\beta^2}$ in \eqref{lem3:6}.
\end{proof}

\textbf{The proof of Lemma \ref{lem4}.}
\begin{proof}{Proof}
		\noindent
By \eqref{update-y}, we immediately obtain
\begin{equation}\label{lem4:2}
\left\| \beta\left( y_k-\operatorname{Prox}_{g,\mathcal{Y}}^{\beta}\left( y_k+ \frac{1}{\beta} \nabla _y\mathcal{L}\left( x_k,y_k,\lambda_k\right) \right) \right) \right\|=\beta\|y_{k+1}-y_k\|.
\end{equation}
On the other hand, by \eqref{update-x} and the Cauchy-Schwartz inequality, we have
\begin{align}\label{lem4:3}
&\|	\alpha ( x_k-\operatorname{Prox}_{h,\mathcal{X}}^{\alpha}( x_k-\frac{1}{\alpha}\nabla _x\mathcal{L}\left( x_k,y_k,\lambda_k \right) ) )\|\nonumber\\
\le &\alpha \| \operatorname{Prox}_{h,\mathcal{X}}^{\alpha}( x_k-\frac{1}{\alpha} \nabla _xL(x_{k},y_{k+1},\lambda_k) )-\operatorname{Prox}_{h,\mathcal{X}}^{\alpha}( x_k-\frac{1}{\alpha} \nabla _xL(x_{k},y_{k},\lambda_k) )\|+\alpha\| x_{k+1}-x_k\|\nonumber\\
\leq  &\alpha\|x_{k+1}-x_k\| +L\|y_{k+1}-y_k\|,
\end{align}
where the last inequality is by the nonexpansive property of the projection operator and Assumption \ref{ass:Lip}. By \eqref{update-lambda}, the Cauchy-Schwartz inequality and the nonexpansive property of the projection operator,  we obtain
\begin{align}\label{lem4:4}
&\|\nabla_{\lambda}\mathcal{L}(x_k,y_k,\lambda_k)\|\nonumber\\
\le &\left\|\nabla_{\lambda}L(x_{k+1},y_{k+1},\lambda_k)-\nabla_{\lambda}L(x_k,y_k,\lambda_k)\right\|+\frac{1}{\gamma}\| \lambda_{k+1}-\lambda_k\|\nonumber\\
\leq &\frac{1}{\gamma}\|\lambda_{k+1}-\lambda_k\|+\|A\|\|x_{k+1}-x_k\| +\|B\|\|y_{k+1}-y_k\|.
\end{align}
Combing \eqref{lem4:2}, \eqref{lem4:3} and \eqref{lem4:4}, and using Cauchy-Schwarz inequality, we complete the proof.
\end{proof}

\textbf{The proof of Theorem \ref{thm1}.}
\begin{proof}{Proof}
		\noindent
By \eqref{thm1:1}, it can be easily checked that $d_1>0$. By multiplying $d_1$ on the both sides of \eqref{lem4:1} and using \eqref{lem3:1} in Lemma \ref{lem3}, we get
\begin{align}
d_1\|\nabla G^{\alpha,\beta,\gamma}( x_k,y_{k},\lambda_k)\|^2
\le S(x_{k},y_{k},\lambda_{k})-S(x_{k+1},y_{k+1},\lambda_{k+1}).\label{thm1:5}
\end{align}
Summing both sides of \eqref{thm1:5} from $k=1$ to $T(\varepsilon)$, we then obtain
\begin{align}
&\sum_{k=1}^{T\left( \varepsilon \right)}{d_1\|\nabla G^{\alpha,\beta,\gamma}( x_k,y_{k},\lambda_k) \| ^2}
\le S(x_{1},y_{1},\lambda_{1})-S(x_{T(\varepsilon)+1},y_{T(\varepsilon)+1},\lambda_{T(\varepsilon)+1}).\label{thm1:6}
\end{align}
Next we proof $\underbar{S}$ is a finite value. 
By the definition of $\Phi(x,\lambda)$ and $S(x,y,\lambda)$ in Lemma \ref{lem3}, we have
$
S(x,y,\lambda)
=2\max\limits_{y\in\mathcal{Y}}\mathcal{L}(x,y,\lambda)-\mathcal{L}(x,y,\lambda)+h(x)+g(y)
\ge\max\limits_{y\in\mathcal{Y}}\mathcal{L}(x,y,\lambda)+h(x)+g(y).
$
Then, we immediately get
\begin{align}
\min_{\lambda}\min_{x\in\mathcal{X}}\min_{y\in\mathcal{Y}}S(x,y,\lambda)
\ge&\min_{\lambda}\min_{x\in\mathcal{X}}\min_{y\in\mathcal{Y}}
\left\{\max_{y\in\mathcal{Y}}\mathcal{L}(x,y,\lambda)+h(x)+g(y)\right\}\nonumber\\
\ge&\min_{\lambda}\min_{x\in\mathcal{X}}\min_{y\in\mathcal{Y}}\left\{\max_{y\in\mathcal{Y}}\mathcal{L}(x,y,\lambda)+h(x)\right\}
+\min_{\lambda}\min_{x\in\mathcal{X}}\min_{y\in\mathcal{Y}}g(y)\nonumber\\
=&\min_{\lambda}\min_{x\in\mathcal{X}}\max_{y\in\mathcal{Y}}\left\{\mathcal{L}(x,y,\lambda)+h(x)\right\}
+\min_{y\in\mathcal{Y}}g(y),\label{thm1:6.5}
\end{align}
where the second inequality is by the fact that $\min_x\left\{f_1(x)+f_2(x)\right\}\ge\min_x f_1(x)+\min_x f_2(x)$. By \eqref{thm1:6.5} and the fact that $\max_y f_1(y)+\max_y f_2(y)\ge\max_y\{f_1(y)+f_2(y)\}$, we get
\begin{align}
\min_{\lambda}\min_{x\in\mathcal{X}}\min_{y\in\mathcal{Y}}S(x,y,\lambda)-2\min_{y\in\mathcal{Y}}g(y)
\ge&\min_{\lambda}\min_{x\in\mathcal{X}}\max_{y\in\mathcal{Y}}\left\{\mathcal{L}(x,y,\lambda)+h(x)-g(y)\right\}\nonumber\\
=&\min_{x\in\mathcal{X}}\max_{\substack{y\in\mathcal{Y}\\Ax+By= c}}\{f(x,y)+h(x)-g(y)\} >-\infty,\label{thm1:7}
\end{align}
where the second last equality is by Theorem \ref{dual}, and the last inequality is by the assumption that  $\mathcal{X}$ and $\mathcal{Y}$ are convex and compact sets. Also, we have $\min_{y\in\mathcal{Y}}g(y)>-\infty$, and thus $\min\limits_{\lambda}\min\limits_{x\in\mathcal{X}}\min\limits_{y\in\mathcal{Y}}S(x,y,\lambda)>-\infty$.  Then, by \eqref{thm1:6} we obtain
$$
\sum_{k=1}^{T\left( \varepsilon \right)}{d_1\|\nabla G^{\alpha,\beta,\gamma}( x_k,y_{k},\lambda_k) \| ^2}\le S(x_{1},y_{1},\lambda_{1})-\underbar{S}=d_2.
$$
In view of the definition of $T(\varepsilon)$, the above inequality implies that $\varepsilon ^2\le \frac{d_2}{T( \varepsilon )d_1}$ or equivalently, $T\left( \varepsilon \right) \le \frac{d_2}{\varepsilon ^2 d_1}$. This completes the proof. 
\end{proof}


\subsection{The Proofs of Lemmas and Theorems in Section 2.4}\label{ap:c}

\textbf{The proof of Lemma \ref{sc:lem1}.}
\begin{proof}{Proof}
		\noindent
Similar to the proof of Lemma \ref{lem3}, by replacing $\alpha$ with $\alpha_k$, $\gamma$ with $\gamma_k$, $\mu$ with $\rho_k$, $\eta$ with $\eta_k$, respectively, we have
\begin{align}
&M_k(x_{k+1},y_{k+1},\lambda_{k+1})-M_k( x_{k},y_{k},\lambda_k)\nonumber\\
\le&-\left(\alpha_k-\frac{L^3}{(L+\beta)^2}-\frac{L(L+\beta)^2\eta_k^2}{\beta^2}
-\frac{L^2}{\rho_k}-\frac{3L}{2}\right)\|x_{k+1}-x_k\|^2 \nonumber \\
&-\left(\frac{1}{\gamma_k}-\frac{2\|B\|^2(L+\beta)^2\eta_k^2}{L\beta^2}-L-\frac{L^2}{\rho_k} \right)\|\lambda_{k+1}-\lambda_{k}\|^2-\left(\beta-\frac{5L}{2}\right)\|y_{k+1}-y_k\|^2.\label{sc:lem1:2}
\end{align}
On the other hand, by \eqref{sc:1} and \eqref{sc:3}, we obtain
\begin{align}
&M_{k+1}(x_{k+1},y_{k+1},\lambda_{k+1})-M_k( x_{k+1},y_{k+1},\lambda_{k+1})\nonumber \\
=&2(\Psi_{k+1}(x_{k+1},\lambda_{k+1})-\Psi_{k}(x_{k+1},\lambda_{k+1}))-L_{k+1}(x_{k+1},y_{k+1},\lambda_{k+1})\nonumber\\
&+L_k(x_{k+1},y_{k+1},\lambda_{k+1}))\nonumber\\
				=&2(\mathcal{L}_{k+1}(x_{k+1},\tilde{y}_{k+1}^*(x_{k+1},\lambda_{k+1}),\lambda_{k+1})-\mathcal{L}_{k}(x_{k+1},\tilde{y}_k^*(x_{k+1},\lambda_{k+1}),\lambda_{k+1}))+\frac{\rho_{k+1}-\rho_k}{2}\|y_{k+1}\|^2\nonumber\\
\le&2(\mathcal{L}_{k+1}(x_{k+1},\tilde{y}_{k+1}^*(x_{k+1},\lambda_{k+1}),\lambda_{k+1})-\mathcal{L}_{k}(x_{k+1},\tilde{y}_{k+1}^*(x_{k+1},\lambda_{k+1}),\lambda_{k+1}))+\frac{\rho_{k+1}-\rho_k}{2}\|y_{k+1}\|^2\nonumber\\
\le&(\rho_k-\rho_{k+1})\sigma_y^2,\label{sc:lem1:3}
\end{align}
where the last inequality holds since $\{\rho_k\}$ is a nonnegative monotonically decreasing sequence. The proof is completed by adding \eqref{sc:lem1:2} and \eqref{sc:lem1:3}.
\end{proof}

\textbf{The proof of Theorem \ref{sc:thm1}.}
\begin{proof}{Proof}
	\noindent
	Firstly, we denote 	
\begin{equation*}
\nabla \bar{G}^{\alpha_k,\beta,\gamma_k}_k\left( x_k,y_k,\lambda_k \right) =\left( \begin{array}{c}
\alpha_k\left( x_k-\operatorname{Prox}_{h,\mathcal{X}}^{\alpha_k}\left( x_k-\frac{1}{\alpha_k}\nabla _x\mathcal{L}\left( x_k,y_k,\lambda_k \right) \right) \right)\\
\beta\left( y_k-\operatorname{Prox}_{g,\mathcal{Y}}^{\beta}\left( y_k+ \frac{1}{\beta} \nabla _y\mathcal{L}_k\left( x_k,y_k,\lambda_k\right) \right) \right)\\
\frac{1}{\gamma_k}\left( \lambda_k-\mathcal{P}_{\Lambda}\left(\lambda_k-\gamma_k\nabla_{\lambda}\mathcal{L}(x_k,y_k,\lambda_k)\right)\right)
\end{array} \right).
\end{equation*}
It then can be easily checked that$\|\nabla G^{\alpha_k,\beta,\gamma_k}( x_{k},y_{k},\lambda_{k}) \|\le\|\nabla \bar{G}^{\alpha_k,\beta,\gamma_k}_k\left( x_{k},y_{k},\lambda_{k} \right) \|+\rho_{k}\|y_{k}\|.$
Next, we give the upper bound of $\|\nabla \bar{G}^{\alpha_k,\beta,\gamma_k}_k\left( x_k,y_k,\lambda_k \right)\|$. Similar to the proof of Lemma \ref{lem4}, by replacing $\alpha$ with $\alpha_k$, $\gamma$ with $\gamma_k$, respectively, we have
\begin{align}
&\|\nabla  \bar{G}^{\alpha_k,\beta,\gamma_k}_k( x_k,y_{k},\lambda_k)\|^2\nonumber\\
\le& (\beta^2+2L^{2}+3\|B\|^2)\|y_{k+1}-y_k\|^2+(2\alpha_k^2+3\|A\|^2)\|x_{k+1}-x_k\|^2+\frac{3}{\gamma_k^2}\|\lambda_{k+1}-\lambda_k\|^2.\label{sc:thm1:3}
\end{align}
Let $\vartheta_k=\frac{L(L+\beta)^2\eta_k^2(\tau-1)}{\beta^2}
+\frac{L^2(\tau-1)}{\rho_k}$ with $\tau>1$. It can be easily checked that $\alpha_k-\frac{L^3}{(L+\beta)^2}-\frac{L(L+\beta)^2\eta_k^2}{\beta^2}
-\frac{L^2}{\rho_k}-\frac{3L}{2}=\vartheta_k$ and $\frac{1}{\gamma_k}-\frac{2\|B\|^2(L+\beta)^2\eta_k^2}{L\beta^2}-L-\frac{L^2}{\rho_k}=\vartheta_k$ by the settings of $\alpha_k$ and $\gamma_k$. In view of Lemma \ref{sc:lem1}, then we immediately obtian
\begin{align}
&\vartheta_k\|x_{k+1}-x_{k}\|^2+\vartheta_k\|\lambda_{k+1}-\lambda_{k}\|^2+\left(\beta-\frac{5L}{2}\right)\| y_{k+1}-y_{k} \|^2\nonumber\\
\le& M_{k}(x_{k},y_{k},\lambda_k)-M_{k+1}(x_{k+1},y_{k+1},\lambda_{k+1})+(\rho_k-\rho_{k+1})\sigma_y^2.\label{sc:thm1:4}
\end{align}
It follows from the definition of $\tilde{d}_1$ that 
\begin{align}
\frac{2\alpha_k^2+3\|A\|^{2}}{\vartheta _k^2}
\leq& \frac{\frac{8L^6}{(L+\beta)^4}+\frac{8L^2(L+\beta)^4\eta_k^4\tau^2}{\beta^4}+\frac{8\tau^2L^4}{\rho_k^2}+12L^2+3\|A\|^2}{\frac{L^2(L+\beta)^4\eta_k^4(\tau-1)^2}{\beta^4}}
\nonumber\\
\leq& \frac{\frac{8L^6}{(L+\beta)^4}+\frac{8L^2(L+\beta)^4\eta_k^4\tau^2}{\beta^4}+\frac{8\tau^2L^4}{\rho_k^2}+12L^2+3\|A\|^2}{\frac{L^2(L+\beta)^4\eta_k^4(\tau-1)^2}{\beta^4}}
\nonumber\\
=& \frac{8\tau^2}{(\tau-1)^2} + \frac{\frac{8L^6\beta^4}{(L+\beta)^4}+12L^2\beta^4+3\|A\|^2\beta^4}{L^2(L+\beta)^4\eta_k^4(\tau-1)^2}+\frac{8\tau^2L^4\beta^4}{L^2(L+\beta)^4\eta_k^4\rho_k^2(\tau-1)^2}
\le\tilde{d}_1,
\end{align}
where the last inequality holds since $\eta_k>\frac{2(L+\beta)}{\rho_k}=k^{1/4}\ge1$. It also follows from the definition of $\tilde{d}_2$ that $\forall k\ge 1$, we have $\frac{3}{\gamma_k^2}\le\tilde{d}_2\vartheta_k$. 
Setting $d_k^{(2)}=\frac{1}{\max\{\frac{\beta^2+2L^{2}+3\|B\|^2}{\beta-\frac{5L}{2}}, \max\{\tilde{d}_1, \tilde{d}_2\}\vartheta_k\}}$, then multiplying $d_k^{(2)}$ on the both sides of \eqref{sc:thm1:3} and combining \eqref{sc:thm1:4}, we then conclude that
\begin{align}
&d_k^{(2)}\|\nabla  \bar{G}^{\alpha_k,\beta,\gamma_k}_k( x_k,y_{k},\lambda_k)\|^2
\le M_{k}(x_{k},y_{k},\lambda_k)-M_{k+1}(x_{k+1},y_{k+1},\lambda_{k+1})+(\rho_k-\rho_{k+1})\sigma_y^2.\label{sc:thm1:6}
\end{align}
Denoting $\tilde{T}_1(\varepsilon)=\min\{k \mid \| \nabla \bar{G}^{\alpha_k,\beta,\gamma_k}_k(x_k,y_{k},\lambda_k) \| \leq \frac{\varepsilon}{2}, k\geq 1\}$. By summing both sides of \eqref{sc:thm1:6} from $k=1$ to $\tilde{T}_1(\varepsilon)$, we obtain
\begin{align}
&\sum_{k=1}^{\tilde{T}_1(\varepsilon)}d_k^{(2)}\|\nabla  \bar{G}^{\alpha_k,\beta,\gamma_k}_k( x_k,y_{k},\lambda_k)\|^2\nonumber\\
\le& M_{1}(x_{1},y_{1},\lambda_1)-M_{\tilde{T}_1(\varepsilon)+1}(x_{\tilde{T}_1(\varepsilon)+1},y_{\tilde{T}_1(\varepsilon)+1},\lambda_{\tilde{T}_1(\varepsilon)+1})+\rho_1\sigma_y^2.\label{sc:thm1:7}
\end{align}
Next, we prove $\underline{M}$ is a finite value. By the definition of $M_k(x,y,\lambda)$ in Lemma \ref{sc:lem1} and  the definition of $\Psi_k(x,\lambda)$, we have
\begin{align}
M_k(x,y,\lambda)+\max_{y\in \mathcal{Y}}\frac{\rho_k}{2}\|y\|^2
\ge&\max_{y\in \mathcal{Y}}\mathcal{L}_k(x,y,\lambda)+\max_{y\in \mathcal{Y}}\frac{\rho_k}{2}\|y\|^2+h(x)+g(y)\nonumber\\
\ge&\max_{y\in \mathcal{Y}}\mathcal{L}(x,y,\lambda)+h(x)+g(y),\label{sc:thm1:8}
\end{align}
where the last inequality is by the fact that $\max_y f_1(y)+\max_y f_2(y)\ge\max_y\{f_1(y)+f_2(y)\}$. Simlilar to the proof of Theorem \ref{thm1}, we conclude $\min\limits_{\lambda}\min\limits_{x\in\mathcal{X}}\min\limits_{y\in\mathcal{Y}}M_k(x,y,\lambda)>-\infty$.  We then conclude from \eqref{sc:thm1:7} that
\begin{align}
\sum_{k=1}^{\tilde{T}_1(\varepsilon)}d_k^{(2)}\|\nabla  \bar{G}^{\alpha_k,\beta,\gamma_k}_k( x_k,y_{k},\lambda_k)\|^2
\le M_{1}(x_{1},y_{1},\lambda_1)-\underline{M}+\rho_1\sigma_y^2=\tilde{d}_3.\label{sc:thm1:9}
\end{align}
Note that $\vartheta_k\ge\frac{L^2(\tau-1)}{\rho_1}$, then we have $\tilde{d}_4\ge\max\left\{\frac{\beta^2+2L^{2}+3\|B\|^2}{(\beta-\frac{5L}{2})\vartheta_k}, \max\{\tilde{d}_1, \tilde{d}_2\}\right\}$ by the definition of $\tilde{d}_4$ which implies that $d_k^{(2)}\ge\frac{1}{\tilde{d}_4\vartheta_k}$. By multiplying $\tilde{d}_4$ on the both sides of \eqref{sc:thm1:9}, we have
$$\frac{\varepsilon^2}{4}\le\frac{\tilde{d}_3\tilde{d}_4}{\sum_{k=1}^{\tilde{T}_1(\varepsilon)}\frac{1}{\vartheta_k}}.$$
Note that $\rho_k=\frac{2(L+\beta)}{k^{1/4}}\le2(L+\beta)$, and thus by the definition of $\eta_k$, we have $\eta_k\le\frac{2(L+2\beta)(L+\beta)}{\rho_k\beta}$. Then we obtain $\sum_{k=1}^{\tilde{T}_1(\varepsilon)}\frac{1}{\vartheta_k}\ge\sum_{k=1}^{\tilde{T}_1(\varepsilon)}\frac{2\beta^4(L+\beta)k^{-1/2}}{L(\tau-1)[2(L+\beta)^3(L+2\beta)^2+L\beta^4]}$. Using the fact that $\sum_{k=1}^{\tilde{T}_1( \varepsilon)}1/k^{1/2}\ge\tilde{T}_1( \varepsilon)^{1/2}$, we conclude that
$$
\tilde{T}_1( \varepsilon)\le \left(\frac{2\tilde{d}_3\tilde{d}_4L(\tau-1)[2(L+\beta)^3(L+2\beta)^2+L\beta^4]}{\beta^4(L+\beta)\varepsilon^2}\right)^{2}.
$$
On the other hand, if $k\ge\frac{256(L+\beta)^4\sigma_y^4}{\varepsilon^4}$, then $\rho_k\le\frac{\varepsilon}{2\sigma_y}$. This inequality together with the definition of $\sigma_y$ then imply that $\rho_k\|y_k\|\le\frac{\varepsilon}{2}$. Therefore, there exists a
$$
T( \varepsilon)\le \max\left\{\frac{4\tilde{d}_3^2\tilde{d}_4^2L^2(\tau-1)^2[2(L+\beta)^3(L+2\beta)^2+L\beta^4]^2}{\beta^8(L+\beta)^2\varepsilon^4},\frac{256(L+\beta)^4\sigma_y^4}{\varepsilon^4}\right\}
$$
such that $\|\nabla G^{\alpha_k,\beta,\gamma_k}( x_{k},y_{k},\lambda_{k}) \|\le\|\nabla \bar{G}^{\alpha_k,\beta,\gamma_k}_k\left( x_{k},y_{k},\lambda_{k} \right) \|+\rho_{k}\|y_{k}\|\le\varepsilon$. 
The proof is completed.
\end{proof}

\end{document}